\documentclass[final,12pt]{elsarticle}

\usepackage{subfigure}
\usepackage{mathrsfs}
\usepackage{eucal}
\usepackage{amsmath}
\usepackage{amssymb}
\usepackage{amsfonts}
\usepackage{multirow}
\usepackage{lscape}

\newcommand{\Eref}[1]{Equation (\ref{#1})}
\newcommand{\fref}[1]{Figure (\ref{#1})}

\newcommand{\BB}{\mathbf{B}}
\newcommand{\DD}{\mathbf{D}}

\newcommand{\aaa}{\mathbf{a}}

\newcommand{\bn}{\mathbf{N}}
\newcommand{\bm}{\mathbf{M}}
\newcommand{\bveps}{\boldsymbol{\varepsilon}}

\journal{Computers and Structures}

\begin{document}

\begin{frontmatter}

\title{Linear free flexural vibration of cracked functionally graded plates in thermal environment}

\author[a]{S~Natarajan}
\author[b]{P~M~Baiz}
\author[c]{M~Ganapathi}
\author[c]{P~Kerfriden}
\author[d]{S~Bordas\fnref{label2}\corref{cor1}}

\address [a]{Graduate Student, Cardiff School of Engineering Theoretical, Applied and Computational Mechanics, Cardiff University, Wales, U.K.}
\address[b]{Lecturer, Department of Aeronautics, Imperial College, London, U.K.}
\address[c]{Head-Stress \& DTA, Aerospace Engineering, Mahindra System, Bangalore, India.}
\address[c]{Lecturer, Cardiff School of Engineering Theoretical, Applied and Computational Mechanics, Cardiff University, Wales, U.K.}
\address[d]{
Professor, Cardiff School of Engineering Theoretical, Applied and Computational Mechanics, Cardiff University, Wales, U.K.}

 \fntext[3]{School of Engineering, Theoretical, Applied and Computational Mechanics, Queen's Building, Room S1.03, Cardiff University, CF24 3AA, Wales , U.K. stephane.bordas@alumni.northwestern.edu. Tel. +44 (0)29 20875941. http://www.engin.cf.ac.uk/whoswho/profile.asp?RecordNo=679, http://www.researcherid.com/rid/A-1858-2009}

\begin{abstract}
In this paper, the linear free flexural vibrations of functionally graded material plates with a through center crack is studied using an 8-noded shear flexible element. The material properties are assumed to be temperature dependent and graded in the thickness direction. The effective material properties are estimated using the Mori-Tanaka homogenization scheme. The formulation is developed based on first-order shear deformation theory. The shear correction factors are evaluated employing the energy equivalence principle. The variation of the plates natural frequency is studied considering various parameters such as the crack length, plate aspect ratio, skew angle, temperature, thickness and boundary conditions. The results obtained here reveal that the natural frequency of the plate decreases with increase in temperature gradient, crack length and gradient index.
\end{abstract}


\end{frontmatter}

\section{Introduction}
A functionally graded material (FGM) is a new class of material whose properties are characterized by the volume fraction of its constituent materials. The concept of characterization of material properties is not new, but the unique feature of FGMs is that these materials are made of a mixture of ceramics and metals that are characterized by the \emph{smooth and continuous} variation in properties from one surface to another~\cite{Koizumi1993,Koizumi1997,Suresh1997}. For structural integrity, FGMs are preferred over fiber-matrix composites that may result in debonding due to the mismatch in the mechanical properties across the interface of two discrete materials bonded together. With the increased use of these materials for structural components in many engineering applications, it is necessary to understand the dynamic characteristics of functionally graded plates. 

Many researchers have recently attempted to study the bending behavior of FGM plates on three-dimensional elasticity solutions~\cite{Reddy2001,Vel2002,Kashtalyan2004,Pan2003}. All these works are limited to simply supported plates under sinusoidal transverse mechanical or thermal loading. Reddy and Cheng~\cite{Reddy2001} and Vel and Batra~\cite{Vel2002} have accounted for the variation of material properties through the thickness according to a power-law distribution and the locally effective material properties were obtained in terms of the volume fractions of the constituents through the Mori-Tanaka homogenization scheme. Kashtalayan~\cite{Kashtalyan2004} derived the elasticity solutions making use of the Plevako general solution of the equilibrium equations for inhomogeneous isotropic media, whereas, Pan~\cite{Pan2003} studied the laminated functionally graded simply supported rectangular plates under sinusoidal surface load, extending the Pagano's solutions which may not be valid for finding the solutions of such plate problems with continuous inhomogeneity. Elishakoff and Gentilini~\cite{Elishakoff2005} investigated the three-dimensional static analysis of clamped functionally graded plates under uniformly distributed load applying the Ritz energy method.

Application of 3D analysis, in general, is quite cumbersome while dealing with complex loading and boundary conditions. Hence, the analysis of isotropic, composite and FGM plates is carried out numerically as well as analytically using plate theories assuming plane stress conditions. Such approximation can predict global displacement and bending moments with sufficient accuracy~\cite{Rohwer1998,Birman2007}. Few analytical and finite element studies on the bending analysis of FGM plates are recently available in the literature using plate theories. Qian~\cite{Qian2004} and Matsunaga~\cite{Matsunaga2009} examined the bending of thick square FGM plates considering the higher-order shear deformation theory, whereas, Zenkour~\cite{Zenkour2006,Zenkour2007} dealt with 2D trigonometric functions based shear deformation theory. The nonlinear thermo-mechanical response of FGM plates was examined by Praveen and Reddy~\cite{Praveen1998} and Reddy~\cite{Reddy2000} considering the higher-order structural theory. Carrera~\cite{Carrera2008} have obtained closed form and finite element solutions for the static analysis of functionally graded plates subjected to transverse mechanical loads. The unified formulation employed in~\cite{Carrera2008} permits a large variety of plate models with variable kinematic assumptions covering first-order as well as higher-order theories. However, higher-order models that involve additional displacement fields may be based on either an equivalent single layer theory or discrete layer approach. Furthermore, they are computationally expensive in the sense that the number of unknowns to be solved is high compared to that of the first-order shear deformation formulation.

In literature, there has been many studies on dynamic characteristics of FGM plates~\cite{Praveen1998,Yang2001,Yang2002,Vel2004,Sundararajan2005,Sundararajan2005a,Ganapathi2006a,Lee2010} and shells~\cite{Sundararajan2007,Sundararajan2007b,Sundararajan2007a,Natarajan2006,Prakash2007}.
The dynamic characteristics of a cracked structural element is especially important because a crack in a vibrating structure results in stiffness decrease, stress concentration, anisotropy and local flexibility, which are functions of the crack location and size. Moreover the crack will open and close depending on the vibration amplitude. The vibration of cracked isotropic plates was studied as early as 1969 by Lynn and Kumbasar~\cite{Lynn1967} who used a Green's function approach. Later, in 1972, Stahl and Keer~\cite{Stahl1972} studied the vibration of cracked rectangular plates using elasticity methods. The other numerical methods that are used to study the dynamic response and instability of plates with cracks or local defects are:  (1) Finite fourier series transform~\cite{Solecki1985}; (2) Rayleigh-Ritz Method~\cite{Khadem2000}; (3) harmonic balance method~\cite{Wu2005}; (4) finite element method~\cite{Qian1991,Lee1993}; (5) extended finite element method ~\cite{Belytschko2009}; (6) Smoothed finite element methods~\cite{Liu2007,Liu2008,Nguyen-Xuan2008b,Liu2009a,Bordas2009a} and (7) Meshfree methods~\cite{Nguyen2008,Qian2004,Ferreira2006,Gilhooley2007}. Recently, Yang and Chen~\cite{Yang2008} and Kitipornchai \textit{et al.,}~\cite{Kitipornchai2009} studied the dynamic characteristics of FGM beams with an edge crack. In case of FGM, the dynamic characterisitcs also depend on the gradient index compared to the isotropic case. However, to the author's knowledge, studies of cracked FGM plates are scarce in the literature that are of practical importance to the designers.

The main focus of this paper is to compute the linear free flexural vibrations of FGM plates with a through center crack using the finite element method.
Here, an eight-noded shear flexible quadrilateral plate element based on field consistency approach~\cite{Prathap1988,Ganapathi1991} is used to study the dynamic characteristics of FGM plates subjected to thermo-mechanical loadings. 
The shear correction factor is calculated from the energy equivalence principle.
The temperature field is assumed to be constant in the plate and varied only in the thickness direction. The material is assumed to be temperature dependent and graded in the thickness direction according to a power law distribution in terms of the volume fractions of the constituents. The effective material properties are estimated from the volume fractions and the material properties of the constituents using the Mori-Tanaka homogenization method~\cite{Mori1973,Benvensite1987}. The formulation developed herein is validated considering different problems for which the solutions are available in the literature. Detailed parametric studies are carried out to understand the free vibration characteristics of cracked FGM plates.

The paper is organized as follows. In Section \ref{FGMTheory}, a brief introduction to FGM is given followed by the plate formulation. The treatment of boundary conditions for skew plates and eight-noded quadrilateral plate element are discussed in Section \ref{eledes}. A detailed numerical study is presented in Section \ref{numeexamples} followed by conclusion in the last section.

\section{Theoretical development and formulation}
\label{FGMTheory}

\subsection{Functionally Graded Material}
A functionally graded material (FGM) rectangular plate (length $a$, width $b$ and thickness $h$), made by mixing two distinct material phases: a metal and ceramic is considered with coordinates $x,y$ along the in-plane directions and $z$ along the thickness direction (see \fref{fig:platefig}). The material on the top surface $(z=h/2)$ of the plate is ceramic and is graded to metal at the bottom surface of the plate $(z=-h/2)$ by a power law distribution. The homogenized material properties are computed using the Mori-Tanaka Scheme~\cite{Mori1973,Benvensite1987}. 

\subsubsection{Estimation of mechanical and thermal properties}
Based on the Mori-Tanaka homogenization method, the effective bulk modulus $K$ and shear modulus $G$ of the FGM are evaluated as~\cite{Mori1973,Benvensite1987,Cheng2000,Qian2004}

\begin{eqnarray}
{K - K_m \over K_c - K_m} &=& {V_c \over 1+(1-V_c){3(K_c - K_m) \over 3K_m + 4G_m}} \nonumber \\
{G - G_m \over G_c - G_m} &=& {V_c \over 1+(1-V_c){(G_c - G_m) \over G_m + f_1}}
\label{eqn:bulkshearmodulus}
\end{eqnarray}

\noindent where

\begin{equation}
f_1 = {G_m (9K_m + 8G_m) \over 6(K_m + 2G_m)}
\end{equation}

\noindent Here, $V_i~(i=c,m)$ is the volume fraction of the phase material. The subscripts $c$ and $m$ refer to the ceramic and metal phases, respectively. The volume fractions of the ceramic and metal phases are related by $V_c + V_m = 1$, and $V_c$ is expressed as

\begin{equation}
V_c(z) = \left( {2z + h \over 2h} \right)^n, \hspace{0.2cm}  n \ge 0
\label{eqn:volFrac}
\end{equation}

\noindent where $n$ in \Eref{eqn:volFrac} is the volume fraction exponent, also referred to as the gradient index. The effective Young's modulus $E$ and Poisson's ratio $\nu$ can be computed from the following expressions:

\begin{eqnarray}
E = {9KG \over 3K+G} \nonumber \\
\nu = {3K - 2G \over 2(3K+G)}
\label{eqn:young}
\end{eqnarray}

\noindent The effective heat conductivity coefficient $\kappa$ and coefficient of thermal expansion $\alpha$ is given by~\cite{Hatta1985,Rosen1970}

\begin{eqnarray}
{\kappa - \kappa_m \over \kappa_c - \kappa_m} &=& {V_c \over 1+(1-V_c){(\kappa_c - \kappa_m) \over 3\kappa_m}} \nonumber \\
{\alpha - \alpha_m \over \alpha_c - \alpha_m} &=&{  \left( {1 \over K} - {1 \over K_m} \right) \over \left( {1 \over K_c} - {1 \over K_m} \right) }
\label{eqn:heatproperties}
\end{eqnarray}

\noindent The effective mass density $\rho$ is given by the rule of mixtures as~\cite{Vel2004}

\begin{equation}
\rho = \rho_c V_c + \rho_m V_m
\label{eqn:mdensity}
\end{equation}

\subsubsection{Temperature distribution through the thickness}
The material properties
$P$ that are temperature dependent can be written as~\cite{Reddy1998}

\begin{equation}
P = P_o(P_{-1}T^{-1} + 1 + P_1 T + P_2 T^2 + P_3 T^3),
\end{equation}

\noindent where $P_o,P_{-1},P_1,P_2,P_3$ are the coefficients of temperature $T$ and are unique to each constituent material phase. The temperature variation is assumed to occur in the thickness direction only and the temperature field is considered to be constant in the $xy$-plane. In such a case, the temperature distribution along the thickness can be obtained by solving a steady state heat transfer equation

\begin{equation}
-{d \over dz} \left[ \kappa(z) {dT \over dz} \right] = 0, \hspace{0.5cm} T = T_c ~\textup{at}~ z = h/2;~~ T = T_m ~\textup{at} ~z = -h/2
\label{eqn:heat}
\end{equation}

\noindent The solution of \Eref{eqn:heat} is obtained by means of a polynomial series~\cite{Wu2004} as

\begin{equation}
T(z) = T_m + (T_c - T_m) \eta(z,h)
\label{eqn:tempsolu}
\end{equation}

\noindent where

\begin{equation}
\begin{split}
\eta(z,h) = {1 \over C} \left[ \left( {2z + h \over 2h} \right) - {\kappa_{cm} \over (n+1)\kappa_m} \left({2z + h \over 2h} \right)^{n+1} + \right. \\ 
\left. {\kappa_{cm} ^2 \over (2n+1)\kappa_m ^2 } \left({2z + h \over 2h} \right)^{2n+1}
-{\kappa_{cm} ^3 \over (3n+1)\kappa_m ^3 } \left({2z + h \over 2h} \right)^{3n+1} \right. \\ + 
\left. {\kappa_{cm} ^4 \over (4n+1)\kappa_m^4 } \left({2z + h \over 2h} \right)^{4n+1} 
- {\kappa_{cm} ^5 \over (5n+1)\kappa_m ^5 } \left({2z + h \over 2h} \right)^{5n+1} \right] ;
\end{split}
\end{equation}

\begin{equation}
\begin{split}
C = 1 - {\kappa_{cm} \over (n+1)\kappa_m} + {\kappa_{cm} ^2 \over (2 n+1)\kappa_m ^2} 
- {\kappa_{cm} ^3 \over (3n+1)\kappa_m ^3} \\ + {\kappa_{cm} ^4 \over (4n+1)\kappa_m ^4}
- {\kappa_{cm} ^5\over (5n+1)\kappa_m ^5}
\end{split}
\end{equation}

\noindent where $\kappa_{cm} = \kappa_c - \kappa_m$ and $T_c,~T_m$ denote the temperature of the ceramic and metal phases, respectively.


\subsection{Plate formulation}
\begin{figure}
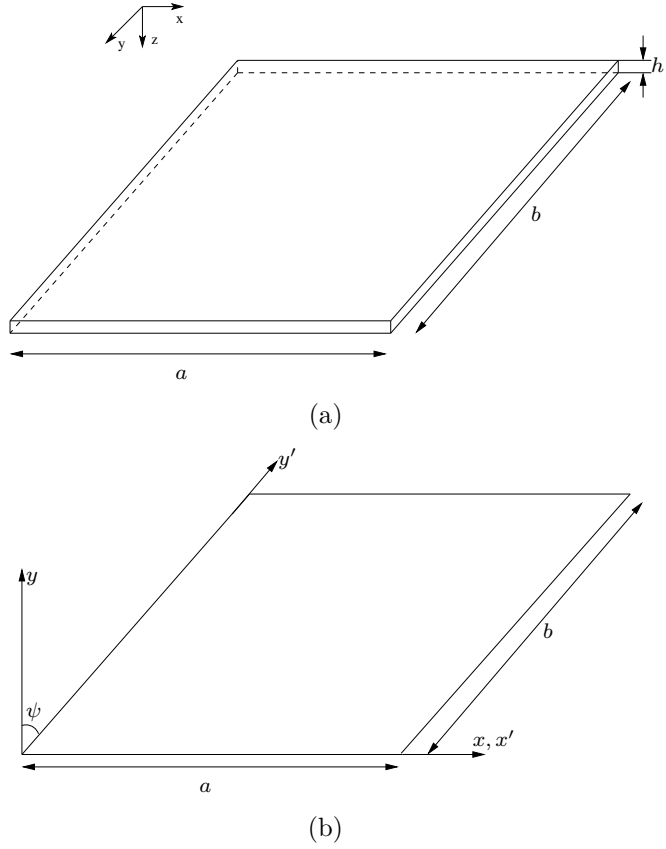

\centering
\subfigure[]{\input{plate.pstex_t}}
\subfigure[]{\input{skew.pstex_t}}
\caption{(a) coordinate system of a rectangular FGM plate, (b) Coordinate system of a skew plate}
\label{fig:platefig}
\end{figure}

Using the Mindlin formulation, the displacements $u,v,w$ at a point $(x,y,z)$ in the plate (see \fref{fig:platefig}) from the medium surface are expressed as functions of the mid-plane displacements $u_o,v_o,w_o$ and indepedent rotations $\theta_x,\theta_y$ of the normal in $yz$ and $xz$ planes, respectively, as

\begin{eqnarray}
u(x,y,z,t) &=& u_o(x,y,t) + z \theta_x(x,y,t) \nonumber \\
v(x,y,z,t) &=& v_o(x,y,t) + z \theta_y(x,y,t) \nonumber \\
w(x,y,z,t) &=& w_o(x,y,t) 
\label{eqn:displacements}
\end{eqnarray}

\noindent where $t$ is the time. The strains in terms of mid-plane deformation can be written as

\begin{equation}
\bveps  = \left\{ \begin{array}{c} \bveps_p \\ 0 \end{array} \right \}  + \left\{ \begin{array}{c} z \bveps_b \\ \bveps_s \end{array} \right\} 
\label{eqn:strain1}
\end{equation}

\noindent The midplane strains $\bveps_p$, bending strain $\bveps_b$, shear strain $\varepsilon_s$ in \Eref{eqn:strain1} are written as

\begin{eqnarray}
\renewcommand{\arraystretch}{1.5}
\bveps_p = \left\{ \begin{array}{c} u_{o,x} \\ v_{o,y} \\ u_{o,y}+v_{o,x} \end{array} \right\}, \hspace{1cm}
\renewcommand{\arraystretch}{1.5}
\bveps_b = \left\{ \begin{array}{c} \theta_{x,x} \\ \theta_{y,y} \\ \theta_{x,y}+\theta_{y,x} \end{array} \right\} \nonumber \\
\renewcommand{\arraystretch}{1.5}
\bveps_s = \left\{ \begin{array}{c} \theta _x + w_{o,x} \\ \theta _y + w_{o,y} \end{array} \right\}, \hspace{1cm}
\renewcommand{\arraystretch}{1.5}
\end{eqnarray}

\noindent where the subscript `comma' represents the partial derivative with respect to the spatial coordinate succeeding it. The membrane stress resultants $\bn$ and the bending stress resultants $\bm$ can be related to the membrane strains, $\bveps_p$ and bending strains $\bveps_b$ through the following constitutive relations

\begin{eqnarray}
\bn &=& \left\{ \begin{array}{c} N_{xx} \\ N_{yy} \\ N_{xy} \end{array} \right\} = \mathbf{A} \bveps_p + \BB \bveps_b - \bn^\textup{T} \nonumber \\
\bm &=& \left\{ \begin{array}{c} M_{xx} \\ M_{yy} \\ M_{xy} \end{array} \right\} = \BB \bveps_p + \DD \bveps_b - \bm^\textup{T} 
\end{eqnarray}

\noindent where the matrices $\aaa = A_{ij}, \BB= B_{ij}$ and $\DD = D_{ij}; (i,j=1,2,6)$ are the extensional, bending-extensional coupling and bending stiffness coefficients and are defined as

\begin{equation}
\left\{ A_{ij}, ~B_{ij}, ~ D_{ij} \right\} = \int_{-h/2}^{h/2} \overline{Q}_{ij} \left\{1,~z,~z^2 \right\}~dz
\end{equation}

\noindent The thermal stress resultant, $N^{\textup{T}}$ and the moment resultant $M^{\textup{T}}$ are

\begin{eqnarray}
\bn^{\textup{T}} &=& \left\{ \begin{array}{c} N_{xx}^{\textup{T}} \\ N_{yy}^{\textup{T}} \\ N_{xy}^{\textup{T}} \end{array} \right\} = \int_{-h/2}^{h/2} \overline{Q}_{ij} \alpha(z,T) \left\{ \begin{array}{c} 1 \\ 1 \\ 0 \end{array} \right\} ~ \Delta T(z) ~ dz \nonumber \\
\bm^{\textup{T}} &=& \left\{ \begin{array}{c} M_{xx}^{\textup{T}} \\ M_{yy}^{\textup{T}} \\ M_{xy}^{\textup{T}} \end{array} \right\} = \int_{-h/2}^{h/2} \overline{Q}_{ij} \alpha(z,T) \left\{ \begin{array}{c} 1 \\ 1 \\ 0 \end{array} \right\} ~z~\Delta T(z) ~ dz
\end{eqnarray}

\noindent where the thermal coefficient of expansion $\alpha(z,T)$ is given by \Eref{eqn:heatproperties} and $\Delta T(z) = T(z)-T_o$ is the temperature rise from the reference temperature $T_o$ at which there are no thermal strains.

Similarly, the transverse shear force $Q = \{Q_{xz},Q_{yz}\}$ is related to the transverse shear strains $\varepsilon_s$ through the following equation

\begin{equation}
Q_{ij} = E_{ij} \varepsilon_s
\end{equation}

\noindent where $E_{ij} = \int_{-h/2}^{h/2} \overline{Q} \upsilon_i \upsilon_j~dz;~ (i,j=4,5)$ is the transverse shear stiffness coefficient, $\upsilon_i, \upsilon_j$ is the transverse shear coefficient for non-uniform shear strain distribution through the plate thickness. The stiffness coefficients $\overline{Q}_{ij}$ are defined as

\begin{eqnarray}
\overline{Q}_{11} = \overline{Q}_{22} = {E(z,T) \over 1-\nu^2}; \hspace{1cm} \overline{Q}_{12} = {\nu E(z,T) \over 1-\nu^2}; \hspace{1cm} \overline{Q}_{16} = \overline{Q}_{26} = 0 \nonumber \\
\overline{Q}_{44} = \overline{Q}_{55} = \overline{Q}_{66} = {E(z,T) \over 2(1+\nu) }
\end{eqnarray}

\noindent where the modulus of elasticity $E(z,T)$ and Poisson's ratio $\nu$ are given by \Eref{eqn:young}. The strain energy function $U$ is given by

\begin{equation}
\begin{split}
U(\boldsymbol{\delta}) = {1 \over 2} \int_{\Omega} \left\{ \bveps_p^{\prime} \mathbf{A} \bveps_p + \bveps_p^{\prime} \mathbf{B} \bveps_b + 
\bveps_b^{\textup{T}} \mathbf{B} \bveps_p + \bveps_b^{\prime} \mathbf{D} \bveps_b +  \bveps_s^{\prime} \mathbf{E} \bveps_s - \right. \\ \left.  \bveps_p^{\textup{T}} \bn - \bveps_b^{\textup{T}} \bm \right\}~ d\Omega
\end{split}
\label{eqn:potential}
\end{equation}

\noindent where $\boldsymbol{\delta} = \{u,v,w,\theta_x,\theta_y\}$ is the vector of the degree of freedom associated to the displacement field in a finite element discretization. Following the procedure given in~\cite{Rajasekaran1973}, the strain energy function $U$ given in~\Eref{eqn:potential} can be rewritten as

\begin{equation}
U(\boldsymbol{\delta}) = {1 \over 2}  \boldsymbol{\delta}^{\textup{T}} \mathbf{K}  \boldsymbol{\delta}
\label{eqn:poten}
\end{equation}

\noindent where $\mathbf{K}$ is the linear stiffness matrix. The kinetic energy of the plate is given by

\begin{equation}
T(\boldsymbol{\delta}) = {1 \over 2} \int_{\Omega} \left\{p (\dot{u}_o^2 + \dot{v}_o^2 + \dot{w}_o^2) + I(\dot{\theta}_x^2 + \dot{\theta}_y^2) \right\}~d\Omega
\label{eqn:kinetic}
\end{equation}

\noindent where $p = \int_{-h/2}^{h/2} \rho(z)~dz, ~ I = \int_{-h/2}^{h/2} z^2 \rho(z)~dz$ and $\rho(z)$ is the mass density that varies through the thickness of the plate given by~\Eref{eqn:mdensity}. The plate is subjected to a temperature field and this in turn results in in-plane stress resultants $(N_{xx}^{\textup{th}}, N_{yy}^{\textup{th}}, N_{xy}^{\textup{th}})$. The external work due to the in-plane stress resultants $(N_{xx}^{\textup{th}}, N_{yy}^{\textup{th}}, N_{xy}^{\textup{th}})$ developed in the plate under the thermal load is

\begin{equation}
\begin{split}
V(\boldsymbol{\delta}) = \int_{\Omega} \left\{ {1 \over 2} \left[ N_{xx}^{\textup{th}} \left( {\partial w \over \partial x} \right)^2 + N_{yy}^{\textup{th}} \left( {\partial w \over \partial y} \right)^2 +
2 N_{xy}^{\textup{th}} \left( {\partial w \over \partial x} \right) \left( {\partial w \over \partial y} \right) \right] \right. \\
\left. + {h^2 \over 24} \left[ N_{xx}^{\textup{th}} \left\{ \left( {\partial \theta _x \over \partial x} \right)^2 + \left( {\partial \theta _y \over \partial x} \right)^2 \right\} + N_{yy}^{\textup{th}} \left\{ \left( {\partial \theta _x \over \partial y} \right)^2 + \left( {\partial \theta _y \over \partial y} \right)^2 \right\} \right. \right. \\
\left. \left. 2N_{xy}^{\textup{th}} \left( {\partial \theta _x \over \partial x} {\partial \theta _x \over \partial y} + {\partial \theta _y \over \partial x} {\partial \theta _y \over \partial y} \right) \right] \right\} ~d\Omega
\end{split}
\label{eqn:potthermal}
\end{equation}

\noindent Substituting \Eref{eqn:poten} - (\ref{eqn:potthermal}) in Lagrange's equations of motion, the following governing equation is obtained

\begin{equation}
\mathbf{M} \ddot{\boldsymbol{\delta}} + (\mathbf{K}+\mathbf{K}_{\rm G}) \boldsymbol{\delta} = 0
\label{eqn:govereqn}
\end{equation}

\noindent where $\mathbf{M}$ is the consistent mass matrix. After substituting the characteristic of the time function~\cite{Ganapathi1991} $\ddot{\boldsymbol{\delta}} = -\omega^2 \boldsymbol{\delta}$, the following algebraic equation is obtained

\begin{equation}
\left( (\mathbf{K}+\mathbf{K}_{\rm G})  - \omega^2 \mathbf{M}\right) \boldsymbol{\delta} = 0
\label{eqn:finaldiscre}
\end{equation}

\noindent where $\mathbf{K}$ and $\mathbf{K}_{\rm G}$ are the stiffness matrix and geometric stiffness matrix due to thermal loads, respectively, $\omega$ is the natural frequency. The plate is subjected to a temperature field. So, the first step in the solution process is to compute the in-plane stress resultants due to the temperature field. These will then be used to compute the stiffness matrix of the system and then the frequencies are computed for the system.

\section{Element Description and Treatment of boundary conditions}
\label{eledes}

\subsection{Eight noded shear flexible element}
The plate element employed here is a $\mathbb{C}_0$ continuous shear flexible element with five degrees of freedom $(u_o,v_o,w_o,\theta_x,\theta_y)$ at eight nodes in an 8-noded quadrilateral (QUAD-8) element. If the interpolation functions for QUAD-8 are used directly to interpolate the five variables $(u_o,v_o,w_o,\theta_x,\theta_y)$ in deriving the shear strains and membrane strains, the element will lock and show oscillations in the shear and membrane stresses. The field consistency requires that the transverse shear strains and membrane strains must be interpolated in a consistent manner. Thus, the $\theta_x$ and $\theta_y$ terms in the expressions for shear strain $\varepsilon_s$ have to be consistent with the derivative of the field functions, $w_{o,x}$ and $w_{o,y}$. This is achieved by using field redistributed substitute shape functions to interpolate those specific terms, which must be consistent as described in~\cite{Prathap1988,Ganapathi1991}. This element is free from locking syndrome and has good convergence properties. For complete description of the element, interested readers are referred to the literature~\cite{Prathap1988,Ganapathi1991}, where the element behavior is discussed in great detail. 

\subsection{Skew boundary transformation}
For skew plates supported on two adjacent edges, the edges of the boundary elements may not be parallel to the global axes $(x,y,z)$. In order to specify the boundary conditions at skew edges, it is necessary to use edge displacements $(u_o^\prime,v_o^\prime,w_o^\prime)$, etc., in local coordinates $(x^\prime,y^\prime,z^\prime)$ (see \fref{fig:platefig}).The element matrices corresponding to the skew edges are transformed from global axes to local axes on which the boundary conditions can be conveniently specified. The relation between the global and local degrees of freedom of a particular node can be obtained through simple transformation rules~\cite{Makhecha2001} expressed as

\begin{equation}
\mathbf{d} = \mathbf{L}_g \mathbf{d}^\prime
\end{equation}

\noindent where $d_i,d_i^\prime$ are the generalized displacement vectors in the global and local coordinate system, repectively, of a node $i$ and are defined as

\begin{eqnarray}
\mathbf{d} &=& \left[ \begin{array}{ccccc}u_o & v_o & w_o & \theta_x & \theta_y \end{array} \right] ^\textup{T} \nonumber \\
\mathbf{d}^\prime &=& \left[ \begin{array}{ccccc}u_o^\prime & v_o^\prime & w_o^\prime & \theta_x^\prime & \theta_y^\prime \end{array} \right] ^\textup{T}
\end{eqnarray}

\noindent The nodal transformation matrix for a node $i$, on the skew boundary is given by

\begin{equation}
\renewcommand\arraystretch{1.5}
\mathbf{L}_g = \left[ \begin{array}{ccccc} \cos \psi & \sin \psi & 0 & 0 & 0 \\
-\sin \psi & \cos \psi & 0 & 0 & 0 \\
0 & 0 & 1 & 0 & 0 \\
0 & 0 & 0 & \cos \psi & \sin \psi \\
0 & 0 & 0 & -\sin \psi & \cos \psi \end{array} \right]
\end{equation}

\noindent where $\psi$ is the angle of the plate. The transformation matrix for the nodes that do not lie on the skew edge is an identity matrix. Thus, for the complete element, the element transformation matrix is written as

\begin{equation}
\mathbf{T}_e = diag \left\langle \begin{array}{cccccccc} \mathbf{L}_g & \mathbf{L}_g & \mathbf{L}_g & \mathbf{L}_g & \mathbf{L}_g & \mathbf{L}_g & \mathbf{L}_g & \mathbf{L}_g \end{array} \right\rangle
\end{equation}

\section{Numerical Examples}
\label{numeexamples}
\begin{figure}[htpb]
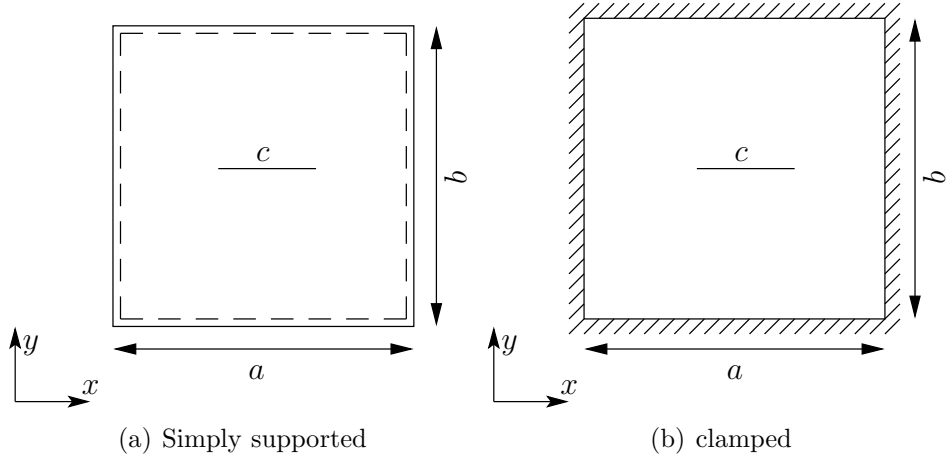

\centering
\subfigure[Simply supported]{\input{simplyS.pstex_t}}
\subfigure[clamped]{\input{clamped.pstex_t}}
\caption{Plate with a center crack}
\label{fig:simplyss}
\end{figure}

In this section, we present the natural frequencies of cracked FGM plates. \fref{fig:simplyss} shows the geometry and boundary condition of the plate. The crack is assumed to be at the center of the plate. For a rectangular plate, the crack is assumed to be parallel to the sides of length $a$. The effect of the gradient index $n$, temperature gradient $(T_c-T_m)$, crack length $c$, aspect ratio $a/b$, thickness of the plate $h$, skew angle of the plate $\psi$  and boundary conditions on the natural frequency are numerically studied. The
assumed values for the parameters are listed in Table \ref{table:parametervalues}.

\begin{table}
\renewcommand\arraystretch{1.5}
\caption{Details of parameters used in the numerical study}
\centering
\begin{tabular}{ll}
\hline
Parameter & Assumed values  \\
\hline
Gradient index, $n$ & 0, 0.5, 1,2, 5, 10 \\
Ceramic Temperature (K), $T_c$ & 300, 400, 600 \\
Metal Temperature (K), $T_m$ & 300, 300, 300 \\
Crack Length $c/a$ & 0.2, 0.4, 0.6 \\
Aspect ratio $a/b$ & 1, 2 \\
Thickness of the plate $a/h$ & 10, 20 \\
Skew angle of the plate $\psi$ (in deg.) & 0$^\circ$, 15$^\circ$, 30$^\circ$, 45$^\circ$ \\
Boundary conditions & Simply supported and Clamped \\
\hline
\end{tabular}
\label{table:parametervalues}
\end{table}

In all cases, we present the non-dimensionalized free flexural frequency defined as

\begin{equation}
\Omega = \omega a^2 \sqrt{ {\rho_c h \over D}}
\end{equation}

\begin{figure}
\includegraphics{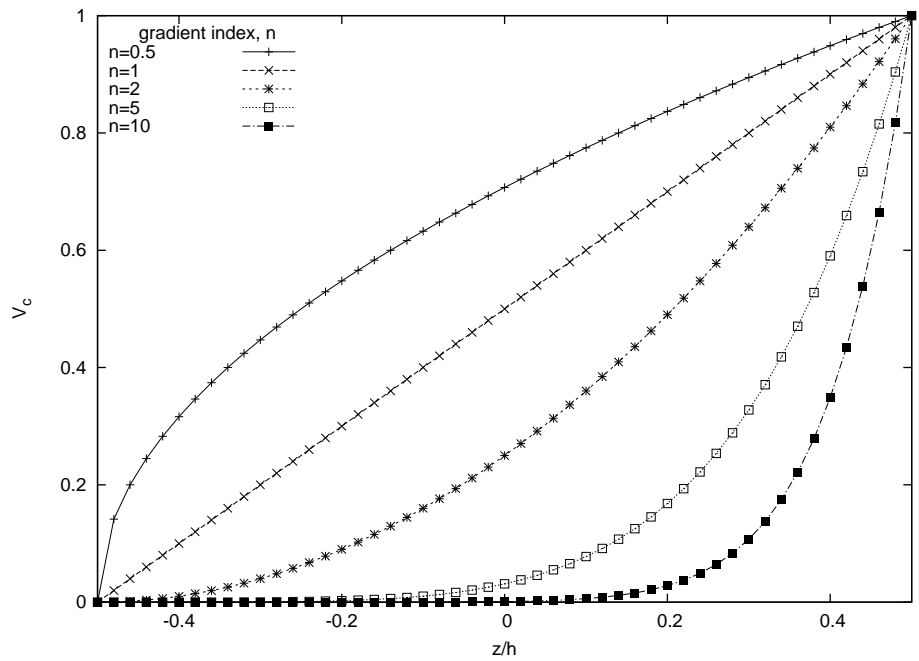}
\caption{Through thickness variation of volume fraction}
\label{fig:volfrac}
\end{figure}

\noindent where $\omega$ is the natural frequency, $D = {E_c h^3 \over 12(1-\nu^2)}$ is the bending rigidity of a homogeneous ceramic plate, $\rho_c,E_c$ are the mass density and Young's modulus of the ceramic, respectively. 
\fref{fig:volfrac} shows the variation of the volume fractions of ceramic and metal, respectively, in the thickness direction $z$ for the FGM plate. The top surface is ceramic rich and the bottom surface is metal rich. 
The FGM plate considered here consists of silicon nitride (SI$_3$N$_4$) and stainless steel (SUS304). The material is considered to be temperature dependent and the temperature coefficients corresponding to SI$_3$N$_4$/SUS304 are listed in Table \ref{table:properties} ~\cite{Sundararajan2005,Reddy1998}. The mass density $(\rho)$ and thermal conductivity $(K)$ are: $\rho_c$=2370 kg/m$^3$, $K_c$=9.19 W/mK for SI$_3$N$_4$ and $\rho_m$ = 8166 kg/m$^3$, $K_m$ = 12.04 W/mK for SUS304. Poisson's ratio $\nu$ is assumed to be constant and taken as 0.28 for the current study~\cite{Sundararajan2005,Sundararajan2005a}. 
Here, the modified shear correction factor obtained based on energy equivalence principle as outlined in~\cite{Singh2011} is used.
The boundary conditions for simply supported and clamped cases are (see \fref{fig:simplyss}):

\noindent \emph{Simply supported boundary condition}: \\
\begin{eqnarray}
u_o = w_o = \theta_y = 0 \hspace{1cm} ~\textup{on} ~ x=0,a \nonumber \\
v_o = w_o = \theta_x = 0 \hspace{1cm} ~\textup{on} ~ y=0,b
\end{eqnarray}

\noindent \emph{Clamped boundary condition}: \\
\begin{eqnarray}
u_o = w_o = \theta_y = v_o = \theta_x = 0 \hspace{1cm} ~\textup{on} ~ x=0,a \nonumber \\
u_o = w_o = \theta_y = v_o = \theta_x = 0 \hspace{1cm} ~\textup{on} ~ y=0,b
\end{eqnarray}

\begin{table}
\renewcommand\arraystretch{1.5}
\caption{Temperature dependent coefficient for material SI$_3$N$_4$/SUS304, Ref. \cite{Reddy1998,Sundararajan2005}.}
\centering
\begin{tabular}{lcccccc}
\hline
Material & Property & $P_o$ & $P_{-1}$ & $P_1$ & $P_2$ & $P_3$  \\
\hline
\multirow{2}{*}{SI$_3$N$_4$} & $E$(Pa) & 348.43e$^9$ &0.0& -3.070e$^{-4}$ & 2.160e$^{-7}$ & -8.946$e^{-11}$  \\
& $\alpha$ (1/K) & 5.8723e$^{-6}$ & 0.0 & 9.095e$^{-4}$ & 0.0 & 0.0 \\
\cline{2-7}
\multirow{2}{*}{SUS304} & $E$(Pa) & 201.04e$^9$ &0.0& 3.079e$^{-4}$ & -6.534e$^{-7}$ & 0.0  \\
& $\alpha$ (1/K) & 12.330e$^{-6}$ & 0.0 & 8.086e$^{-4}$ & 0.0 & 0.0 \\
\hline
\end{tabular}
\label{table:properties}
\end{table}

\noindent Before proceeding with the detailed study on the effect of different parameters on the natural frequency, the formulation developed herein is validated against available results
pertaining to the linear frequencies for functionally graded plate based on higher-order shear deformation plate theory~\cite{Ferreira2006,Qian2004}  and to the linear frequencies for cracked isotropic plates~\cite{Liew1994, Natarajan2010}. The computed frequencies for a square simply supported plate: (a) with gradient index $n=$1 and $a/h=10$ is given in Table~\ref{table:SSValidHigherTheory}; (b) with a center crack and thickness $a/h=1000$ are given in Table \ref{table:SSValidation}.
It can be seen that the numerical results from the present formulation are found to be in good agreement with the existing solutions. 

It is seen that with increase in crack length $c/a$, the natural frequency decreases. The decrease in frequency with increasing crack length $c/a$ is due to reduction in the stiffness of the plate. 
The linear frequencies for uncracked FGM plates in thermal environment is shown in Table~\ref{table:tempValidation} and the results are compared with the results available in the literature~\cite{Huang2004}. Here the calculated non-dimensional frequency is defined as: $\Omega = \omega (a^2/h)\sqrt{ \rho_m(1-\nu^2)/E_m}$. Based on the progressive mesh refinement, an 8$\times$8 mesh is found to be adequate to model the full plate for the uncracked FGM plates. The mesh size used for cracked FGM plates and number of degrees of freedom are given in Table~\ref{table:meshsize}.

\begin{table}[htpb]
\renewcommand\arraystretch{2}
\caption{Non-dimensionalized natural frequency $\Omega = \omega h\sqrt{ \rho_m/E_m}$ for a simply supported square plate for ${a \over h} = 10$, Gradient index, $n$ = 1. The material properties are assumed to be the same as in~\cite{Ferreira2006,Qian2004}. The present formulation based on First order Shear Deformation plate Theory is compared with Third order Shear Deformation Theory (TSDT)~\cite{Ferreira2006} and Higher Order Shear and Normal Deformable Plate Theory (HOSNDPT)~\cite{Qian2004}. } 
\centering
\begin{tabular}{cccc}
\hline Mode & TSDT~\cite{Ferreira2006} & HOSDNPT~\cite{Qian2004} &  FSDT~(Present) \\
  \hline
 1& 0.0592 & 0.0584 & 0.0609 \\
 2&0.1428 & 0.1410 & 0.1422 \\
 3&0.1428 & 0.1410&0.1422\\
\hline
\end{tabular}
\label{table:SSValidHigherTheory}
\end{table}

\begin{table}[htpb]
\renewcommand\arraystretch{2}
\caption{Non-dimensionalized natural frequency for a simply supported square plate with a center crack for ${a \over h} = 1000$. } 
\centering
\begin{tabular}{lrrrrrrr}
\hline c/a & \multicolumn{3}{c}{Mode 1} & &  
\multicolumn{3}{c} {Mode 2} \\
\cline{2-4} \cline{6-8}
  & Ref~\cite{Liew1994} & Ref~\cite{Natarajan2010} & Present &  & Ref~\cite{Liew1994} & Ref~\cite{Natarajan2010} & Present \\
  \hline
 0.0 & 19.74 & 19.74  & 19.752 & &49.35 &  49.38 & 49.410 \\
 0.2 & 19.38 & 19.40 & 19.399 & & 49.16 & 49.85 & 49.300 \\
 0.4 & 18.44 & 18.50 & 18.484 & & 46.44 & 47.27 & 47.018 \\
 0.6 & 17.33 & 17.37 & 17.439 & & 37.75 & 38.92 & 38.928 \\
\hline
\end{tabular}
\label{table:SSValidation}
\end{table}

\begin{table}
\renewcommand\arraystretch{2}
\caption{Comparison of non-dimensional linear frequencies of simply supported FGM plate $(a/b=1, a/h=8)$ in thermal environment.}
\centering
\begin{tabular}{ccrrrrr}
\hline
Temperature & gradient & \multicolumn{2}{c}{Mode 1} & & \multicolumn{2}{c}{Mode 2} \\
\cline{3-4} \cline{5-7}
$T_c,T_m$& $n$ & Ref.~\cite{Huang2004} & Present & & Ref.~\cite{Huang2004} & Present \\
\hline
\multirow{4}{*}{$T_c$=400} & 0.0 & 12.397 & 12.311 & & 29.083 & 29.016 \\
$T_m$=300, & 0.5 & 8.615 & 8.276 && 20.215 & 19.772 \\
& 1.0 & 7.474 & 7.302 & & 17.607 & 17.369 \\
& 2.0 & 6.693 & 6.572 & & 15.762 & 15.599 \\
\cline{2-7}
\multirow{4 }{*}{$T_c$=600} & 0.0 & 11.984 & 11.888 & & 28.504 & 28.421 \\
$T_m$=300, & 0.5 & 8.269 & 7.943 & & 19.784 & 19.327 \\
& 1.0 & 7.171 & 6.989 & & 17.213 & 16.959 \\
& 2.0 & 6.398 & 6.269 & & 15.384 & 15.207 \\
\hline
\end{tabular}
\label{table:tempValidation}
\end{table}

\begin{table}
\caption{Number of nodes and total number of degrees of freedom for FGM plates with and without crack used in the present analysis.}
\centering
\begin{tabular}{lrr}
\hline
c/a & Number  & Total degrees\\
& of nodes & of freedom \\
\hline
0.0 & 225 & 513 \\
0.2 & 184 & 840 \\
0.4 & 264 & 1216 \\
0.6 & 356 & 1604 \\
\hline
\end{tabular}
\label{table:meshsize}
\end{table}

\begin{figure}
\centering
\subfigure[$a/h$=10]{\includegraphics[scale=0.80]{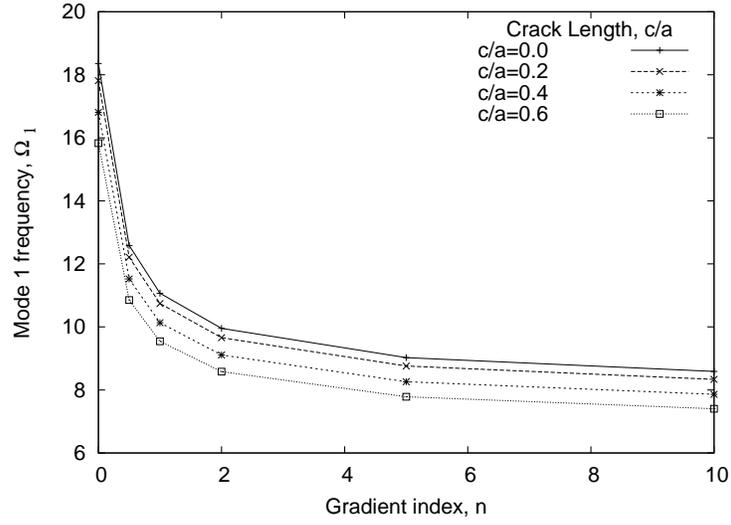}}
\subfigure[$a/h$=20]{\includegraphics[scale=0.80]{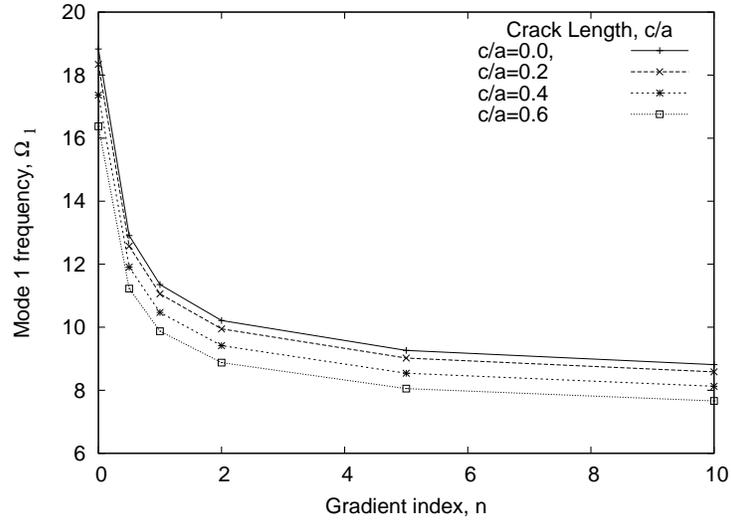}}
\caption{Frequency (Mode 1) as a function of gradient index $n$ and crack length $c/a$ for a simply supported square FGM plate in ambient temperature $(T_c = 300 \textup{K},~ T_m = 300 \textup{K})$.}
\label{fig:modegrph1}
\end{figure}

\begin{figure}
\centering
\subfigure[$a/h$=10]{\includegraphics[scale=0.80]{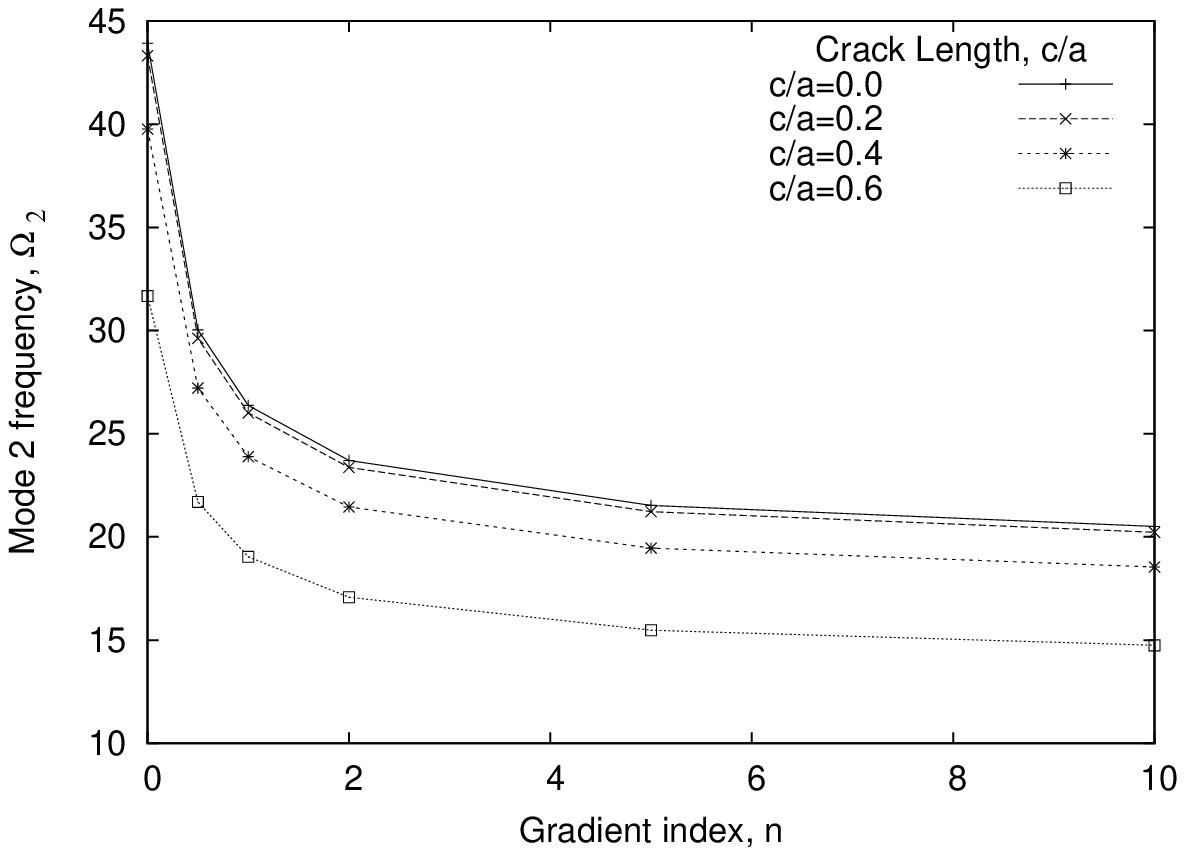}}
\subfigure[$a/h$=20]{\includegraphics[scale=0.80]{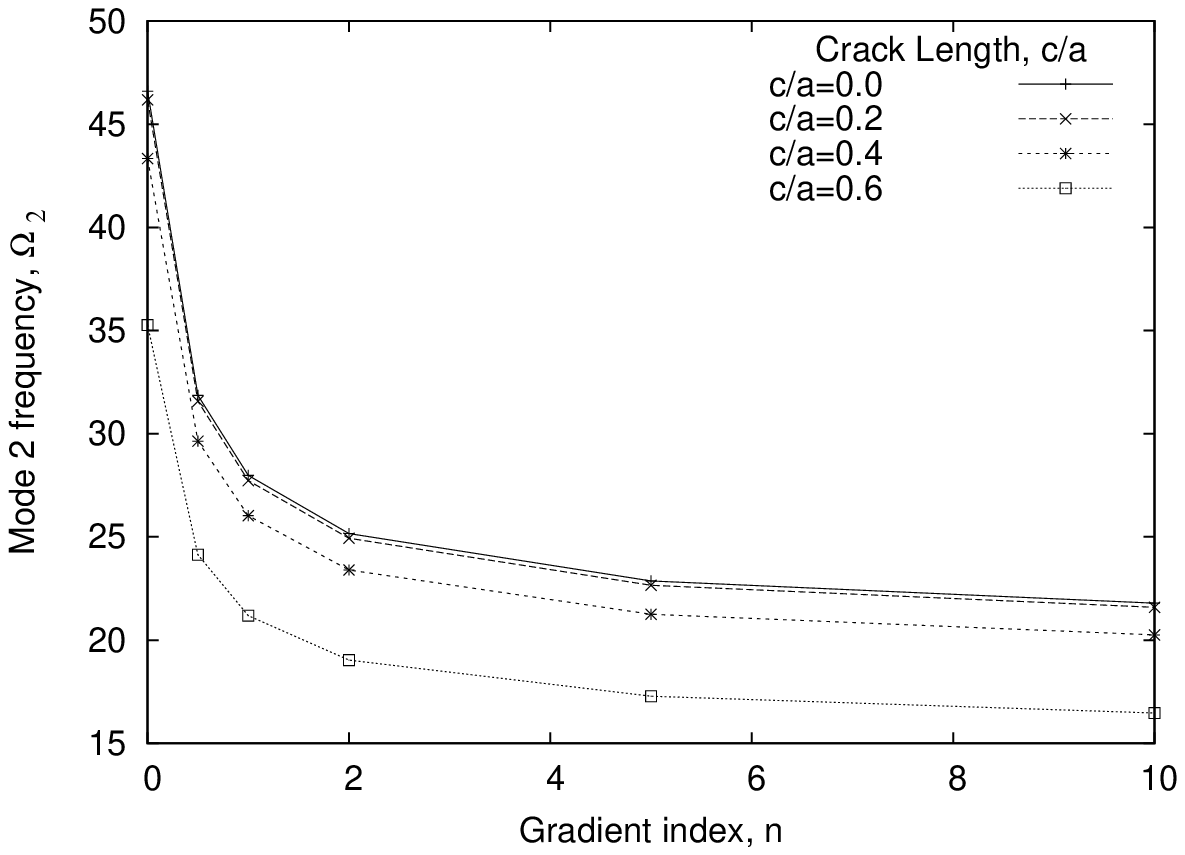}}
\caption{Frequency (Mode 2) as a function of gradient index $n$ and crack length $c/a$ for a simply supported square FGM plate in ambient temperature $(T_c = 300 \textup{K},~ T_m = 300 \textup{K})$.}
\label{fig:modegrph2}
\end{figure}

\begin{figure}
\centering
\subfigure[$a/h$=10]{\includegraphics[scale=0.80]{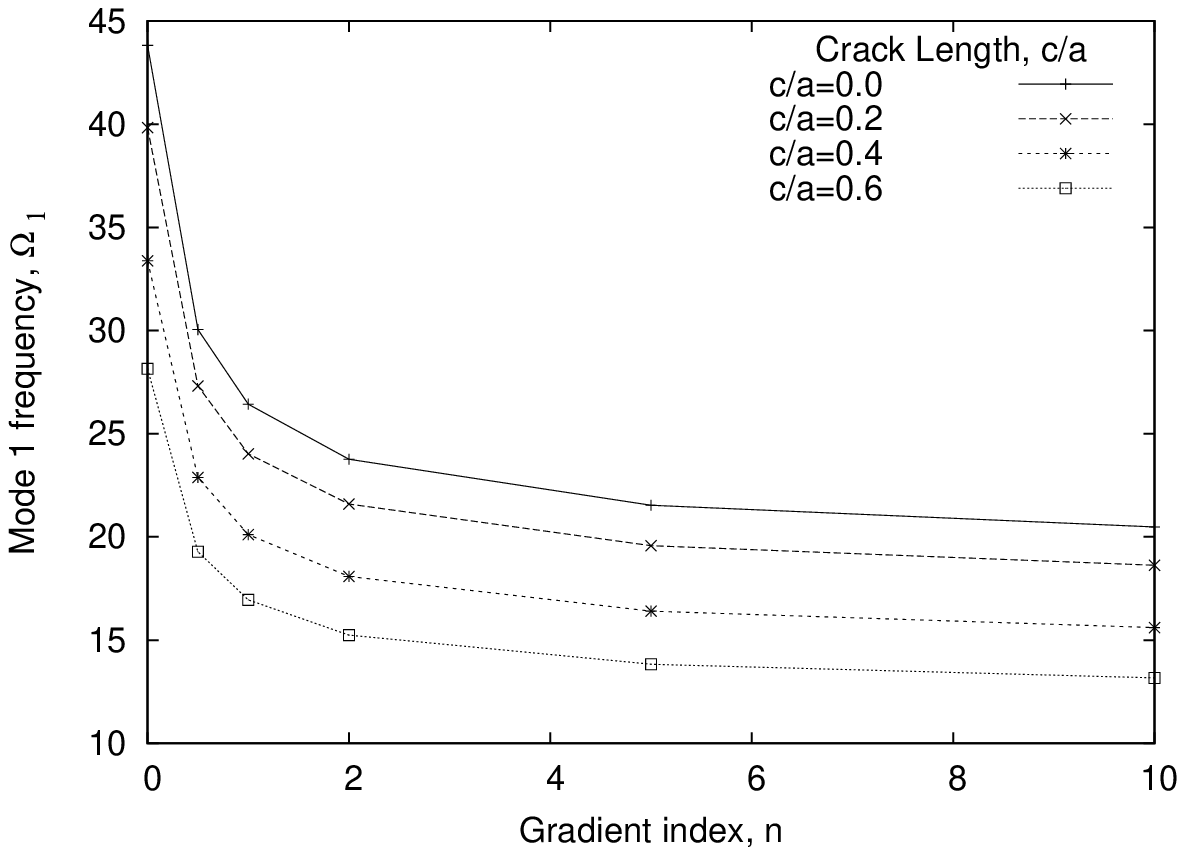}}
\subfigure[$a/h$=20]{\includegraphics[scale=0.80]{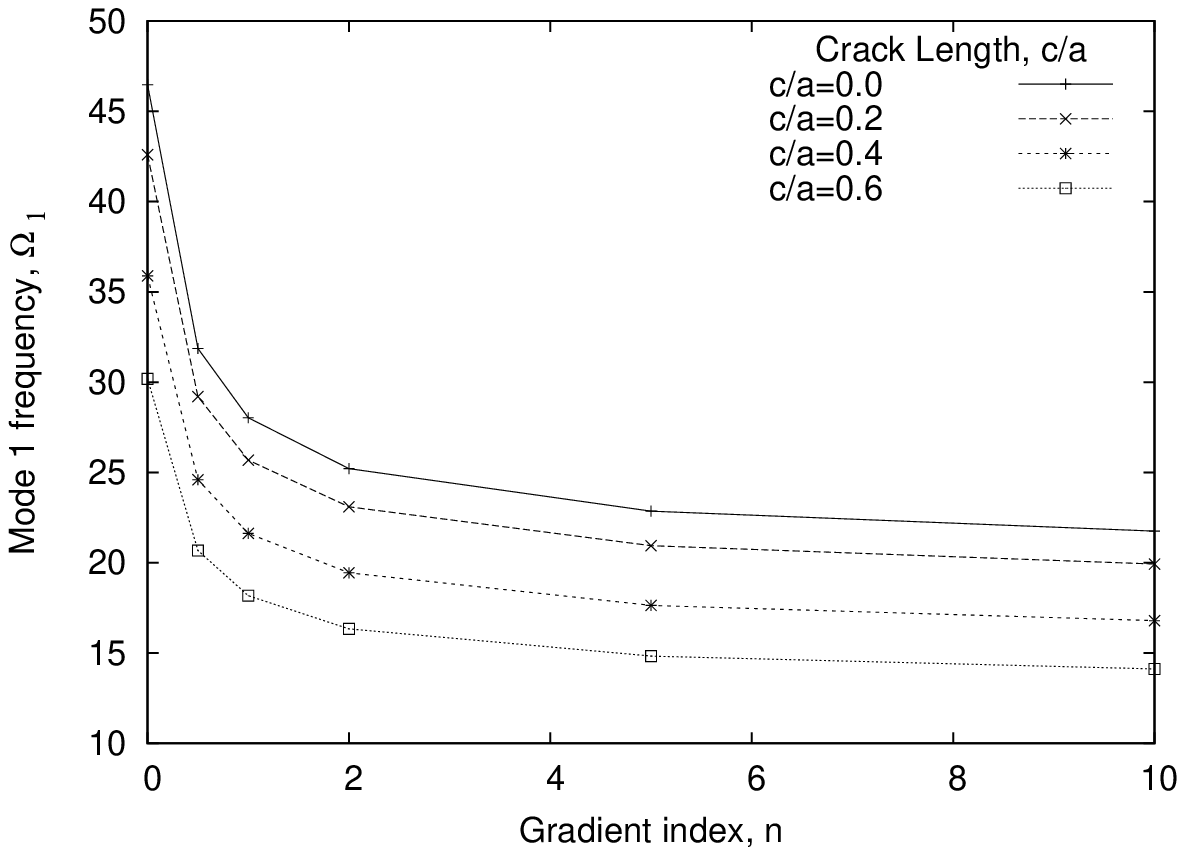}}
\caption{Frequency (Mode 1) as a function of gradient index $n$ and crack length $c/a$ for a simply supported rectangular FGM plate in ambient temperature $(T_c = 300 \textup{K},~ T_m = 300 \textup{K})$.}
\label{fig:modegrph3}
\end{figure}

\begin{figure}
\centering
\subfigure[$a/h$=10]{\includegraphics[scale=0.80]{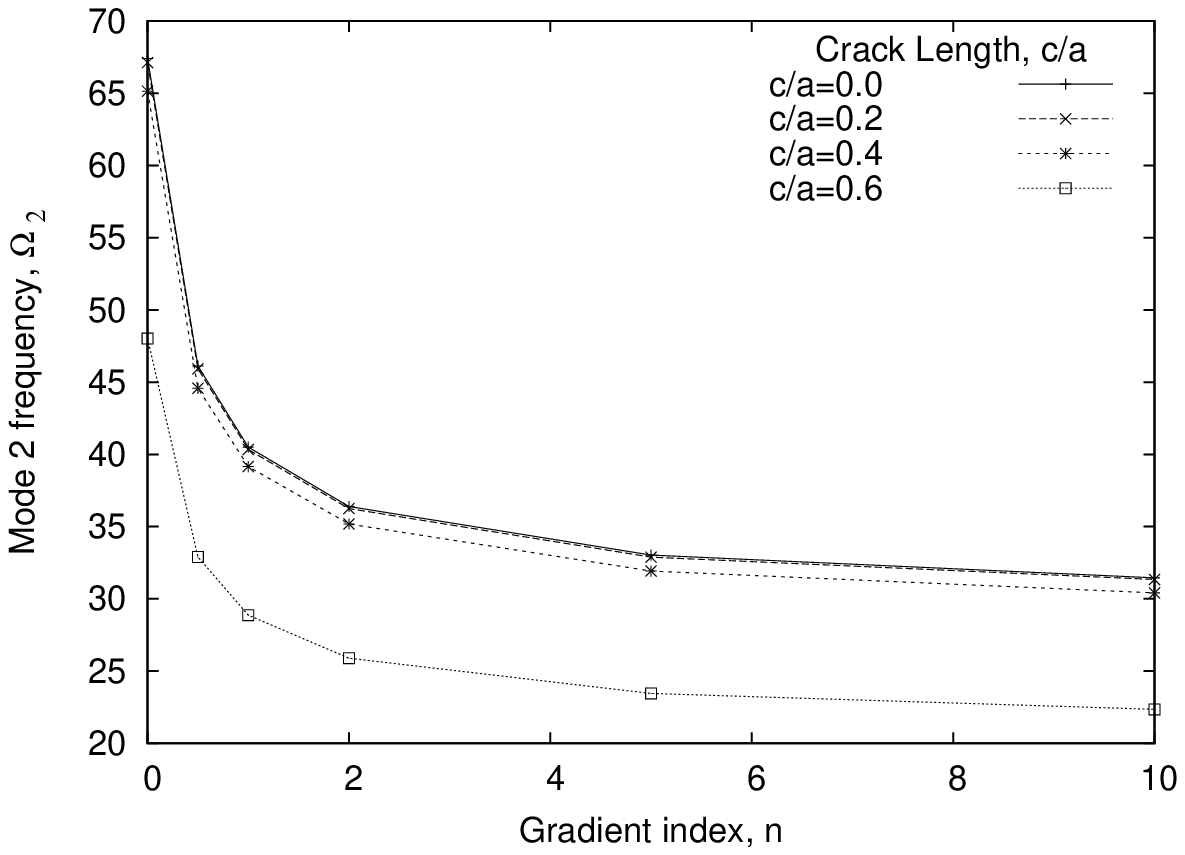}}
\subfigure[$a/h$=20]{\includegraphics[scale=0.80]{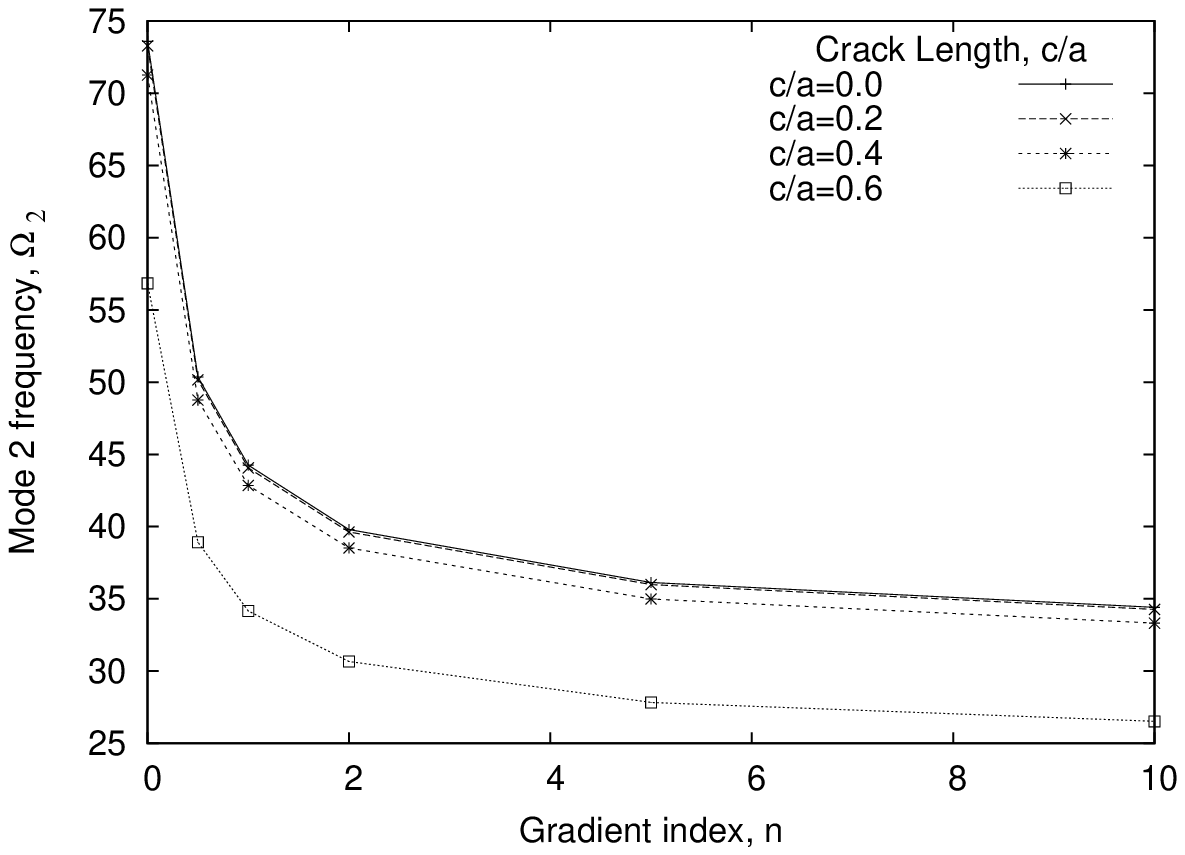}}
\caption{Frequency (Mode 2) as a function of gradient index $n$ and crack length $c/a$ for a simply supported rectangular FGM plate in ambient temperature $(T_c = 300 \textup{K},~ T_m = 300 \textup{K})$.}
\label{fig:modegrph4}
\end{figure}

Next, the linear flexural vibration behavior of FGM is numerically studied with and without thermal environment. For the uniform temperature case, the material properties are evaluated at $T_c=T_m=T=300$K. The variation of natural frequency with aspect ratio, thickness, gradient index and crack length is shown in Table~\ref{table:SSambienTemp} and graphically in Figures (\ref{fig:modegrph1}) - (\ref{fig:modegrph4}). 
\fref{fig:modegrph1} and \fref{fig:modegrph2} shows the variation of mode 1 and mode 2 frequency of a simply supported square plate with gradient index for different crack length $c/a$ and thickness $h$. 
\fref{fig:modegrph3} and \fref{fig:modegrph4} shows the variation of mode 1 and mode 2 frequency of a simply supported rectangular plate with gradient index for different crack length $c/a$ and thickness $h$. The natural frequency increases with decrease in the plate thickness and this behavior is maintained, regardless of the boundary conditions. It is seen from Table~\ref{table:SSambienTemp} and Figures (\ref{fig:modegrph1}) - (\ref{fig:modegrph4}) that with increase in crack length and gradient index, the linear frequency of the plate decreases and with decrease in plate thickness, the natural frequency increases. The decrease in frequency with increasing crack length $c/a$ is due to the reduction in the stiffness of the plate and the decrease in frequency with increasing gradient index $n$ occurs due to increase in the metallic volume fraction. In both the cases, the reduction in frequency is due to the stiffness degradation. It can also be observed that the decrease in frequency is not proportional to the increase in crack length. The frequency decreases rapidly with increase in crack length. This observation is true, irrespective of thick or thin, and square or rectangular case considered in the present study. 

\begin{figure}
\centering
\subfigure[Mode 1]{\includegraphics[scale=0.80]{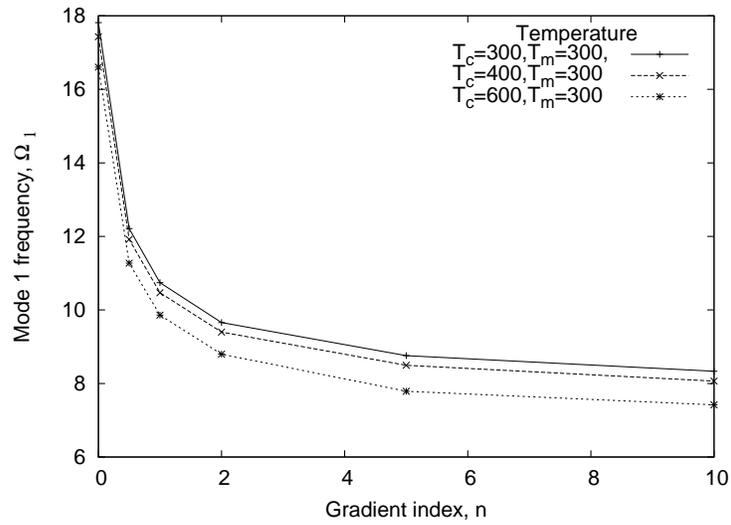}}
\subfigure[Mode 2]{\includegraphics[scale=0.80]{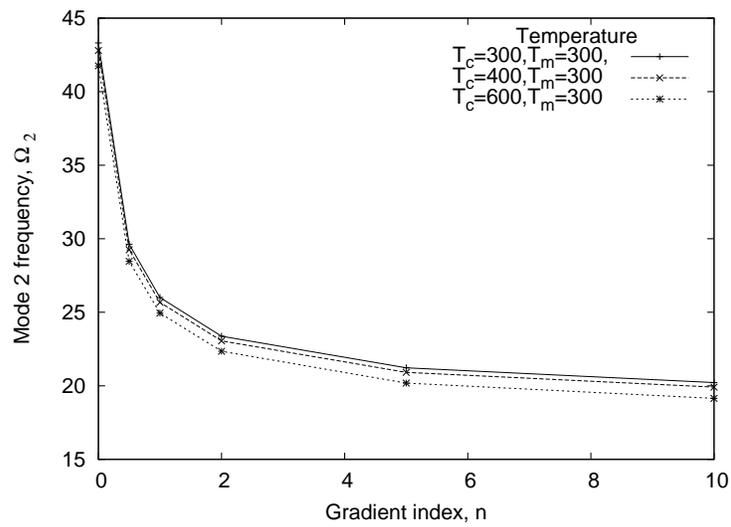}}
\caption{Frequency as a function of gradient index $n$ and different temperature gradient for a simply supported rectangular FGM plate. The results are shown for crack length $c/a=0.2$.}
\label{fig:modegrph5}
\end{figure}

\begin{landscape}
\begin{table}
\renewcommand\arraystretch{1}
\caption{Effect of aspect ratio, thickness, crack length and gradient index on the natural frequency of simply supported FGM plates in ambient temperature $(T_c = 300 \textup{K},~ T_m = 300 \textup{K})$.}
\begin{tabular}{lllrrrrrrrrr}
\hline 
&&& \multicolumn{4}{c}{Mode 1}& & \multicolumn{4}{c}{Mode 2} \\
\cline{4-12}
$a/b$ & $a/h$ & $n$ & \multicolumn{4}{c}{c/a}& & \multicolumn{4}{c}{c/a} \\
\cline{4-7} \cline{9-12}
& & & 0.0 & 0.2 & 0.4 & 0.6 & &0.0 & 0.2 & 0.4 & 0.6 \\
\hline
\multirow{12}{*}{1}&\multirow{6}{*}{10} 
&0.0 & 18.357 & 17.815 & 16.805 & 15.829 & & 43.923 & 43.316 & 39.765 & 31.671 \\
&&0.5 & 12.584 & 12.214 & 11.522 & 10.853 & & 30.032 & 29.622 & 27.218 & 21.697 \\
&&1.0 & 11.069 & 10.742 & 10.132 & 9.544 & & 26.371 & 26.008 & 23.888 & 19.036 \\
&&2.0 & 9.954 & 9.659 & 9.111 & 8.582 & & 23.703 & 23.372 & 21.447 & 17.077 \\
&&5.0 & 9.026  & 8.759  & 8.262 & 7.782 & & 21.527 & 21.222 & 19.455 & 15.479 \\
&&10.0 & 8.588 & 8.334  & 7.861 & 7.405 & & 20.511 & 20.221 & 18.538 & 14.749 \\
\cline{3-12}
&\multirow{6}{*}{20} 
&0.0 & 18.829 & 18.340 & 17.366 & 16.377 & & 46.599 & 46.178 & 43.337 & 35.265 \\
&&0.5 & 12.908 & 12.573 & 11.906 & 11.228 & & 31.862 & 31.575 & 29.643 & 24.133 \\
&&1.0 & 11.353 & 11.058 & 10.470 & 9.873 & & 27.983 & 27.730 & 26.028 & 21.184 \\
&&2.0 & 10.212 & 9.946 & 9.416 & 8.879 & & 25.165 & 24.936 & 23.396 & 19.033 \\
&&5.0 & 9.263  & 9.021  & 8.541 & 8.054 & & 22.868 & 22.658 & 21.252 & 17.281 \\
&&10.0 & 8.813 & 8.584  & 8.127 & 7.664 & & 21.786 & 21.586 & 20.247 & 16.465 \\
\cline{2-12}
\multirow{12}{*}{2}&\multirow{6}{*}{10} 
&0.0 & 43.821 & 39.836 & 33.381 & 28.147 & & 67.409 & 67.127 & 65.149 & 48.012 \\
&&0.5 & 30.052 & 27.315 & 22.877 & 19.282 & & 46.123 & 45.930 & 44.576 & 32.901 \\
&&1.0 & 26.433 & 24.023 & 20.116 & 16.953 & & 40.510 & 40.341 & 39.153 & 28.867 \\
&&2.0 & 23.761 & 21.595 & 18.085 & 15.245 & & 36.397 & 36.246 & 35.181 & 25.888 \\
&&5.0 & 21.532  & 19.572  & 16.399 & 13.830 & & 33.026 & 32.889 & 31.924 & 23.450 \\
&&10.0 & 20.485 & 18.622  & 15.606 & 13.163 & & 31.460 & 31.329 & 30.409 & 22.342 \\
\cline{3-12}
&\multirow{6}{*}{20} 
&0.0 & 46.468 & 42.594 & 35.886 & 30.178 & & 73.591 & 73.283 & 71.260 & 56.844 \\
&&0.5 & 31.864 & 29.205 & 24.596 & 20.676 & & 50.361 & 50.151 & 48.766 & 38.920 \\
&&1.0 & 28.031 & 25.684 & 21.621 & 18.172 & & 44.251 & 44.067 & 42.846 & 34.154 \\
&&2.0 & 25.211 & 23.095 & 19.439 & 16.338 & & 39.792 & 39.626 & 38.526 & 30.667 \\
&&5.0 & 22.860  & 20.943  & 17.634 & 14.826 & & 36.131 & 35.980 & 34.982 & 27.826 \\
&&10.0 & 21.747 & 19.927  & 16.784 & 14.114 & & 34.408 & 34.264 & 33.315 & 26.513 \\
\hline
\end{tabular}
\label{table:SSambienTemp}
\end{table}

\begin{table}
\renewcommand\arraystretch{1}
\caption{Effect of aspect ratio, thickness, crack length and gradient index on the natural frequency of simply supported FGM plates with a temperature gradient $(T_c = 400 \textup{K},~ T_m = 300 \textup{K})$}
\begin{tabular}{lllrrrrrrrrr}
\hline 
&&& \multicolumn{4}{c}{Mode 1}& & \multicolumn{4}{c}{Mode 2} \\
\cline{4-12}
$a/b$ & $a/h$ & $n$ & \multicolumn{4}{c}{c/a}& & \multicolumn{4}{c}{c/a} \\
\cline{4-7} \cline{9-12}
& & & 0.0 & 0.2 & 0.4 & 0.6 & &0.0 & 0.2 & 0.4 & 0.6 \\
\hline
\multirow{12}{*}{1}&\multirow{6}{*}{10} 
&0.0 & 17.962 & 17.433 & 16.465 & 15.537 & & 43.382 & 42.803 & 39.331 & 31.338 \\
&&0.5 & 12.281 & 11.920 & 11.261 & 10.630 & & 29.640 & 29.251 & 26.907 & 21.458 \\
&&1.0 & 10.786 & 10.468 & 9.889 & 9.337 & & 26.014 & 25.672 & 23.606 & 18.820 \\
&&2.0 & 9.862 & 9.397 & 8.877 & 8.383 & & 23.366 & 23.055 & 21.183 & 16.876 \\
&&5.0 & 8.753  & 8.495  & 8.028 & 7.584 & & 21.196 & 20.912 & 19.199 & 15.285 \\
&&10.0 & 8.309 & 8.065  & 7.623 & 7.204 & & 20.177 & 19.909 & 18.282 & 14.556 \\
\cline{3-12}
&\multirow{6}{*}{20} 
&0.0 & 17.533 & 17.095 & 16.275 & 15.463 & & 45.169 & 44.850 & 42.272 & 34.474 \\
&&0.5 & 11.880 & 11.585 & 11.041 & 10.505 & & 30.756 & 30.552 & 28.827 & 23.529 \\
&&1.0 & 10.379 & 10.122 & 9.651 & 9.189 & & 26.947 & 26.772 & 25.267 & 20.622 \\
&&2.0 & 9.255 & 9.026 & 8.612 & 8.209 & & 24.158 & 24.006 & 22.660 & 18.491 \\
&&5.0 & 8.276  & 8.073  & 7.713 & 7.365 & & 21.843 & 21.714 & 20.507 & 16.735 \\
&&10.0 & 7.785 & 7.596  & 7.266 & 6.949 & & 20.730 & 20.614 & 19.483 & 15.905 \\
\cline{2-12}
\multirow{12}{*}{2}&\multirow{6}{*}{10} 
&0.0 & 43.280 & 39.345 & 32.999 & 27.855 & & 66.733 & 66.461 & 64.522 & 47.643 \\
&&0.5 & 29.651 & 26.951 & 22.596 & 19.069 & & 45.642 & 45.456 & 44.132 & 32.647 \\
&&1.0 & 26.067 & 23.690 & 19.860 & 16.760 & & 40.077 & 39.915 & 38.754 & 28.644 \\
&&2.0 & 23.416 & 21.281 & 17.844 & 15.064 & & 35.995 & 35.849 & 34.810 & 25.685 \\
&&5.0 & 21.196  & 19.266  & 16.166 & 13.657 & & 32.638 & 32.507 & 35.168 & 23.262 \\
&&10.0 & 20.149 & 18.317  & 15.374 & 12.992 & & 31.073 & 30.949 & 30.057 & 22.160 \\
\cline{3-12}
&\multirow{6}{*}{20} 
&0.0 & 45.036 & 41.302 & 34.921 & 29.480 & & 72.011 & 71.734 & 69.833 & 56.154 \\
&&0.5 & 30.755 & 28.204 & 23.847 & 20.137 & & 49.155 & 48.968 & 47.677 & 38.412 \\
&&1.0 & 26.992 & 24.746 & 20.921 & 17.667 & & 43.129 & 42.965 & 41.833 & 33.691 \\
&&2.0 & 24.202 & 22.184 & 18.760 & 15.851 & & 38.710 & 38.565 & 37.552 & 30.230 \\
&&5.0 & 21.835  & 20.018  & 16.948 & 14.337 & & 35.044 & 34.914 & 34.006 & 27.400 \\
&&10.0 & 20.691 & 18.974  & 16.079 & 13.614 & & 33.295 & 33.175 & 32.320 & 26.087 \\
\hline
\end{tabular}
\label{table:SSTcer400}
\end{table}

\begin{table}
\renewcommand\arraystretch{1}
\caption{Effect of aspect ratio, thickness, crack length and gradient index on the natural frequency of simply supported FGM plates with a temperature gradient $(T_c = 600 \textup{K},~ T_m = 300 \textup{K})$}
\begin{tabular}{lllrrrrrrrrr}
\hline 
& & & \multicolumn{4}{c}{Mode 1}& & \multicolumn{4}{c}{Mode 2} \\
\cline{4-12}
$a/b$ & $a/h$ & $n$ & \multicolumn{4}{c}{c/a}& & \multicolumn{4}{c}{c/a} \\
\cline{4-7} \cline{9-12}
& & & 0.0 & 0.2 & 0.4 & 0.6 & &0.0 & 0.2 & 0.4 & 0.6 \\
\hline
\multirow{12}{*}{1}&\multirow{6}{*}{10} 
&0.0 & 17.105 & 16.605 & 15.730 & 14.910 & & 42.272 & 41.754 & 38.451 & 30.666 \\
&&0.5 & 11.611 & 11.273 & 10.687 & 10.142 & & 28.808 & 28.467 & 26.252 & 20.958 \\
&&1.0 & 10.155 & 9.859 & 9.350 & 8.878 & & 25.241 & 24.944 & 23.000 & 18.357 \\
&&2.0 & 9.064 & 8.799 & 8.349 & 7.935 & & 22.617 & 22.350 & 20.596 & 16.429 \\
&&5.0 & 8.115  & 7.788  & 7.483 & 7.123 & & 20.425 & 20.187 & 18.598 & 14.828 \\
&&10.0 & 7.642 & 7.420  & 7.054 & 6.723 & & 19.371 & 19.150 & 17.654 & 14.079 \\
\cline{3-12}
&\multirow{6}{*}{20} 
&0.0 & 14.325 & 14.017 & 13.611 & 13.269 & & 41.964 & 41.849 & 39.931 & 32.750 \\
&&0.5 & 9.269 & 9.082 & 8.881 & 8.734 & & 28.239 & 28.165 & 27.002 & 22.190 \\
&&1.0 & 7.866 & 7.712 & 7.573 & 7.492 & & 24.564 & 24.501 & 23.545 & 19.360 \\
&&2.0 & 6.725 & 6.602 & 6.531 & 6.515 & & 21.812 & 21.759 & 20.971 & 17.254 \\
&&5.0 & 5.543  & 5.458  & 5.483 & 5.565 & & 19.401 & 19.358 & 18.755 & 15.455 \\
&&10.0 & 4.818 & 4.761  & 4.864 & 5.023 & & 18.161 & 18.125 & 17.645 & 14.564 \\
\cline{2-12}
\multirow{12}{*}{2}&\multirow{6}{*}{10} 
&0.0 & 42.172 & 38.336 & 32.221 & 27.270 & & 65.381 & 65.129 & 63.276 & 46.934 \\
&&0.5 & 28.805 & 26.181 & 22.005 & 18.626 & & 44.640 & 44.470 & 43.210 & 32.132 \\
&&1.0 & 25.277 & 22.971 & 19.309 & 16.348 & & 39.151 & 39.003 & 37.902 & 28.171 \\
&&2.0 & 22.650 & 20.584 & 17.311 & 14.666 & & 35.101 & 34.970 & 33.989 & 25.234 \\
&&5.0 & 20.412  & 18.553  & 15.621 & 13.252 & & 31.722 & 31.607 & 30.730 & 22.805 \\
&&10.0 & 19.335 & 17.575  & 14.808 & 12.571 & & 30.118 & 30.010 & 29.182 & 21.686 \\
\cline{3-12}
&\multirow{6}{*}{20} 
&0.0 & 41.826 & 38.408 & 32.776 & 27.953 & & 68.585 & 68.377 & 66.757 & 54.742 \\
&&0.5 & 28.234 & 25.928 & 22.160 & 18.933 & & 46.490 & 46.357 & 45.283 & 37.338 \\
&&1.0 & 24.608 & 22.593 & 19.324 & 16.530 & & 50.619 & 40.506 & 39.579 & 32.689 \\
&&2.0 & 21.854 & 20.064 & 17.190 & 14.735 & & 36.253 & 36.157 & 35.346 & 29.260 \\
&&5.0 & 19.387  & 17.809  & 15.318 & 13.185 & & 32.503 & 32.427 & 31.732 & 26.409 \\
&&10.0 & 18.115 & 16.651  & 14.370 & 12.410 & & 30.636 & 30.573 & 29.943 & 25.056 \\
\hline
\end{tabular}
\label{table:SSTcer600}
\end{table}
\end{landscape}

\begin{figure}[htpb]
\centering
\subfigure[Mode 1]{\includegraphics[scale=0.42]{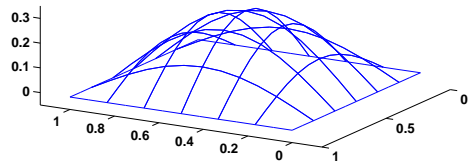}}
\subfigure[Mode 2]{\includegraphics[scale=0.42]{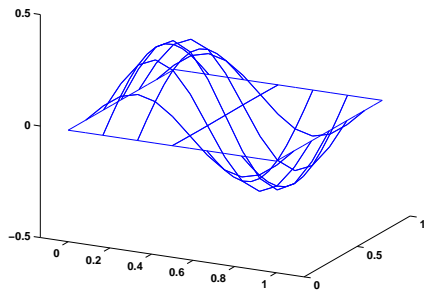}}
\caption{Mode shape for a square plate without a crack. The mode shape is normalized with the maximum value.}
\label{fig:modeshapes}
\end{figure}

\begin{figure}[htpb]
\centering
\subfigure[Mode 1]{\includegraphics[scale=0.42]{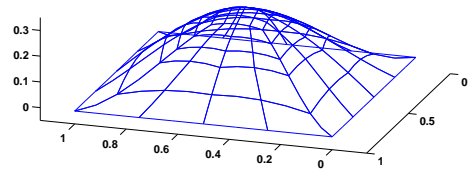}}
\subfigure[Mode 2]{\includegraphics[scale=0.42]{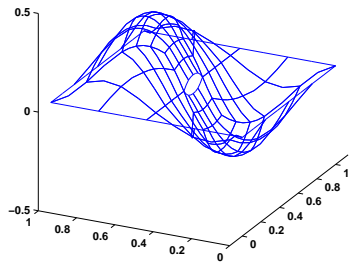}}
\caption{Mode shape for a simply supported square plate with a crack $(c/a=0.2,~a/h=10,~T_c=600,~T_m=300,~n=0)$. The mode shape is normalized with the maximum value.}
\label{fig:crkmodeshapes}
\end{figure}

\begin{figure}
\centering
\subfigure[c/a=0.2]{\includegraphics[scale=0.42]{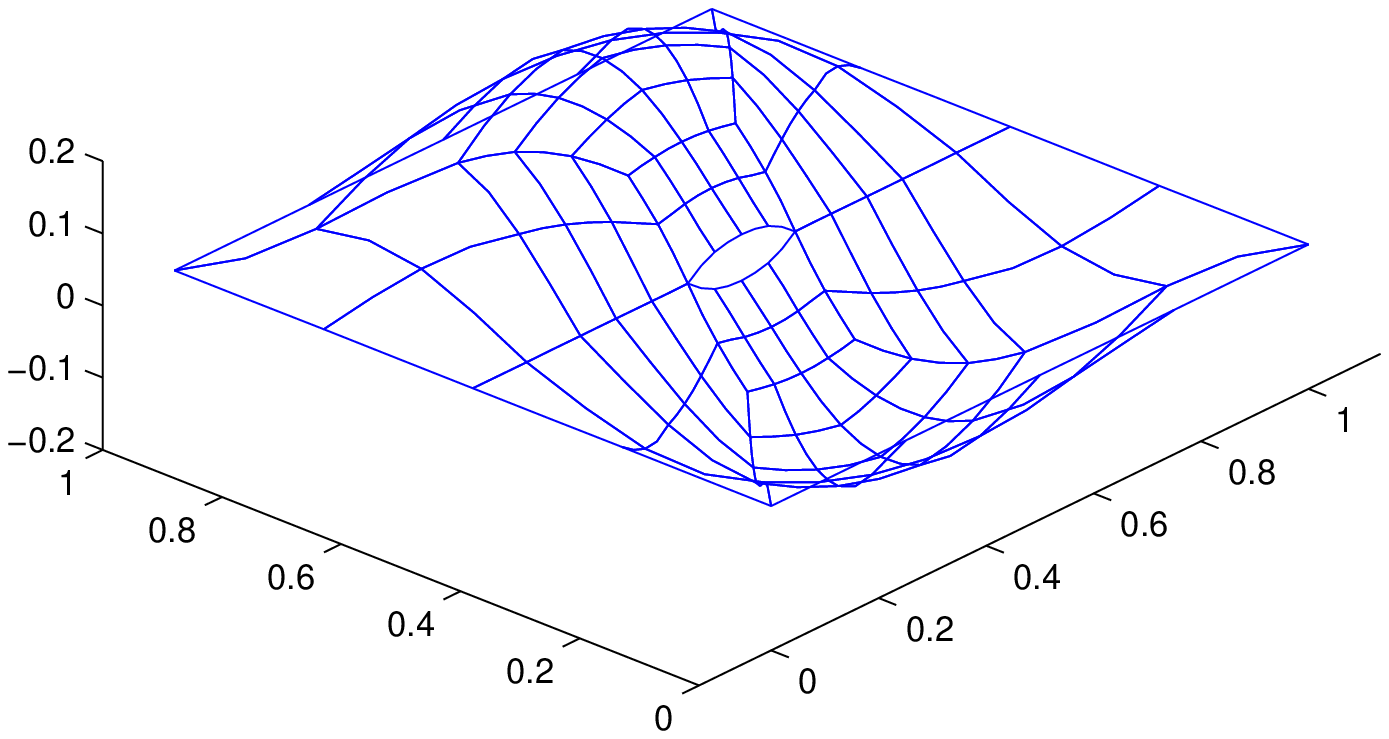}}
\subfigure[c/a=0.4]{\includegraphics[scale=0.42]{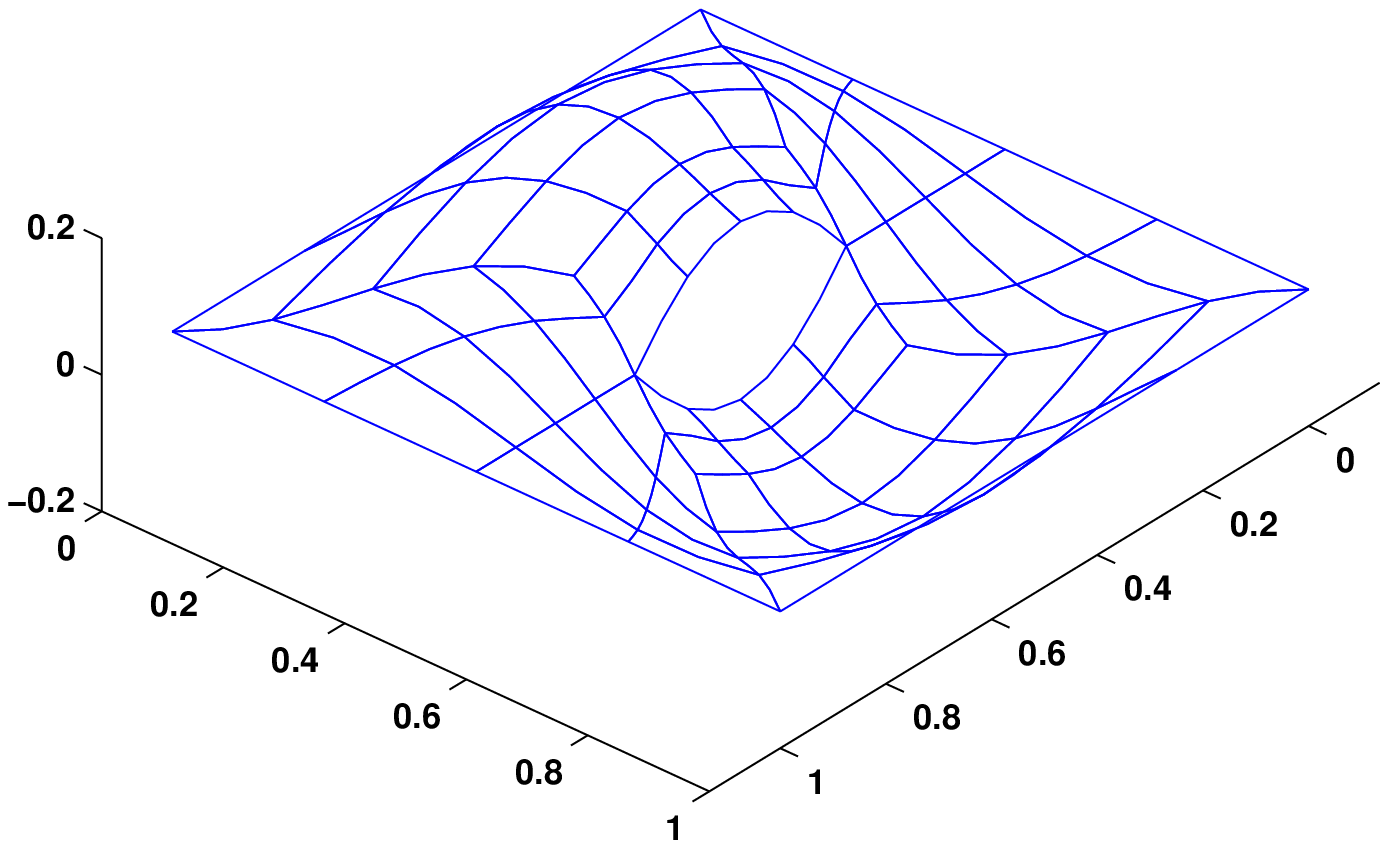}}
\subfigure[c/a=0.6]{\includegraphics[scale=0.45]{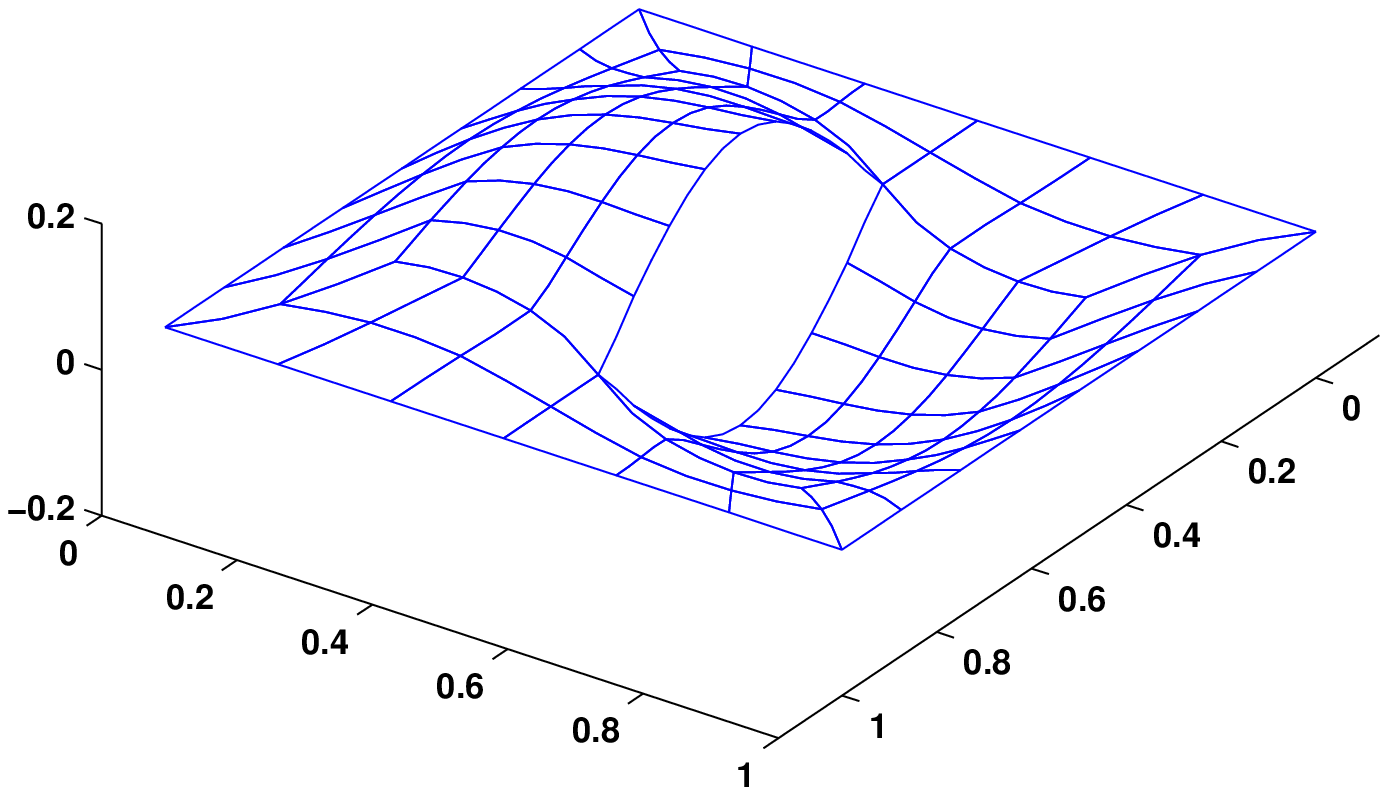}}
\caption{Mode shape (second mode) for a simply supported square plate with different crack sizes $(a/h=10,~T_c=600,~T_m=300,~n=0)$. The mode shape is normalized with the maximum value.}
\label{fig:crkmodeshapes1}
\end{figure}

The influence of thermal gradient on the linear frequencies of FGMs is examined in Tables \ref{table:SSTcer400} and \ref{table:SSTcer600} and graphically in \fref{fig:modegrph5} with different surface temperatures. The temperature field is assumed to vary only in the thickness direction and is determined by~\Eref{eqn:tempsolu}. The temperature for the ceramic surface is varied  $(T_c = 400K, 600K)$ while maintaining a constant value on the metallic surface $(T_m=300K)$. Figures (\ref{fig:modeshapes}) and (\ref{fig:crkmodeshapes}) shows the first two mode shapes of simply supported FGM plates with and without a crack. A similar trend can be observed for the FGM plates in thermal environment. The natural frequency decreases with increasing gradient index and crack length. Also, the frequency further decreases with increase in temperature gradient as expected. The mode shape (second mode) for a simply supported square plate with a center crack is shown in~\fref{fig:crkmodeshapes1} for different crack length and $a/h=10,~T_c=600,~T_m=300,~n=0$. It is seen that with increase in crack length, the plate becomes locally flexible.

\begin{table}
\centering
\renewcommand\arraystretch{1.5}
\caption{Influence of skew angle on linear frequencies for a simply supported FGM plate $(a/b=1, a/h=10)$ with a center crack $(c/a=0.2)$ and with a temperature gradient, $(T_c = 600 \textup{K},~ T_m = 300 \textup{K})$}
\begin{tabular}{clrrrrrr}
\hline 
Skew & Frequency & \multicolumn{6}{c}{gradient index, $n$ } \\
\cline{3-8}
angle, $\psi$ & &0.0 & 0.5 & 1.0 & 2.0 & 5.0 & 10.0 \\
\hline
\multirow{2}{*}{0$^\circ$}  &Mode 1 & 16.605 & 11.273 & 9.859 & 8.799 & 7.788 & 7.420 \\
&Mode 2 &  41.754 & 28.467 & 24.944 & 22.350 & 20.187 & 19.150\\
\cline{3-8}
\multirow{2}{*}{15$^\circ$}   &Mode 1 & 18.351 & 12.480 & 10.918 & 9.751 & 8.748 & 8.257 \\
&Mode 2 &  41.778 & 28.487 & 24.962 & 22.365 & 20.202 & 19.169 \\
\cline{3-8}
\multirow{2}{*}{30$^\circ$}   &Mode 1 &23.956 & 16.341 & 14.310 & 12.796 & 11.519 & 10.911 \\
&Mode 2 &  47.407 & 32.365 & 28.365 & 25.416 & 22.969 & 21.814 \\
\cline{3-8}
\multirow{2}{*}{45$^\circ$}  &Mode 1 &35.793 & 24.473 & 21.448 & 19.200 & 17.325 & 16.450 \\
&Mode 2 &  61.697 & 42.180 & 36.978 & 33.134 & 29.956 & 28.477 \\
\hline
\end{tabular}
\label{table:SSTcer600skew}
\end{table}

\begin{table}
\centering
\renewcommand\arraystretch{1.5}
\caption{Influence of skew angle on linear frequencies for a clamped FGM plate $(a/b=1, a/h=10)$ with a center crack $(c/a=0.2)$ and with a temperature gradient, $(T_c = 600 \textup{K},~ T_m = 300 \textup{K})$}
\begin{tabular}{clrrrrrr}
\hline 
Skew & Frequency & \multicolumn{6}{c}{gradient index, $n$ } \\
\cline{3-8}
angle, $\psi$ & &0.0 & 0.5 & 1.0 & 2.0 & 5.0 & 10.0 \\
\hline
\multirow{2}{*}{0$^\circ$}  &Mode 1 & 29.331 & 19.991 & 17.501 & 15.666 & 14.147 & 13.431 \\
&Mode 2 &  57.569 & 39.137 & 34.283 & 30.691 & 27.726 & 26.335\\
\cline{3-8}
\multirow{2}{*}{15$^\circ$}   & Mode 1 & 30.963 & 21.110 & 18.482 & 16.546 & 14.944 & 14.189 \\
&Mode 2 &  57.338 & 39.178 & 34.318 & 30.724 & 27.755 & 26.382 \\
\cline{3-8}
\multirow{2}{*}{30$^\circ$}   &Mode 1 & 36.649 & 25.010 & 21.903 & 19.611 & 17.717 & 16.831 \\
&Mode 2 &  63.181 & 43.191 & 37.834 & 33.865 & 30.590 & 29.082 \\
\cline{3-8}
\multirow{2}{*}{45$^\circ$}  &Mode 1 & 49.730 & 33.987 & 29.773 & 26.658 & 24.088 & 22.899 \\
&Mode 2 &  78.840 & 53.951 & 106.69 & 42.281 & 38.179 & 36.311 \\
\hline
\end{tabular}
\label{table:ccer600skew}
\end{table}

Lastly, the effect of skewness of a  simply supported and clamped square FGM plates in thermal environment is investigated. The thickness of the plate is taken as $a/h$=10 and the plate is assumed to have a center crack $c/a=0.2$. The results are shown in Tables \ref{table:SSTcer600skew} and \ref{table:ccer600skew} for various gradient indices and skew angles. It can be seen from Tables \ref{table:SSTcer600skew} and \ref{table:ccer600skew} that with increase in skew angle $\psi$, the linear frequency increases and with increase in gradient index $n$, the linear frequency decreases.


\section{Conclusion}
The linear free flexural vibration behavior of FGM plates with and without thermal environment is numerically studied using QUAD-8 shear flexible element. The formulation is based on first order shear deformation theory. The material is assumed to be temperature dependent and graded in the thickness direction according to the power-law distribution in terms of volume fractions of the constituents. The effective material properties are estimated using Mori-Tanaka homogenization method. Numerical experiments have been conducted to bring out the effect of gradient index, aspect ratio, crack length, thickness and boundary condition on the natural frequency of the FGM plate. From the detailed parametric study, it can be concluded that with increase in both the gradient index and crack length, the natural frequency decreases. In both cases, the decrease is due to stiffness degradation. 


\noindent {\bf Acknowledgements} \\

S Natarajan acknowledges the financial support of (1) the
Overseas Research Students Awards Scheme; (2) the Faculty of
Engineering, University of Glasgow, for period Jan. 2008 - Sept. 2009 and of (3) the School
of Engineering (Cardiff University) for the period Sept. 2009 onwards.

S Bordas would like to acknowledge the financial support of the Royal Academy
of Engineering and of the Leverhulme Trust for his Senior Research
Fellowship entitled "Towards the next generation surgical simulators"
as well as the support of EPSRC under grants EP/G069352/1 Advanced
discretisation strategies for "atomistic" nano CMOS simulation and
EP/G042705/1 Increased Reliability for Industrially Relevant Automatic
Crack Growth Simulation with the eXtended Finite Element Method.

\bibliographystyle{elsarticle-num}
\bibliography{myRefFGM}

\end{document}